\documentclass[11pt]{amsart}   	
\usepackage[english]{babel}
\usepackage{geometry}                		
\usepackage{graphicx}				
\usepackage{amsthm}								
\usepackage{amssymb}
\usepackage{amsmath}
\usepackage{tikz-cd}
\usepackage{tikz}
\usetikzlibrary{calc}
\usepackage{hyperref}
\usepackage{multirow}
\usepackage{booktabs}
\usepackage[foot]{amsaddr}

\usetikzlibrary{positioning}

\DeclareMathOperator{\Cat}{Cat}
\DeclareMathOperator{\Kirk}{Kirk}
\DeclareMathOperator{\Fix}{Fix}

\DeclareMathOperator{\NC}{NC}
\DeclareMathOperator{\rk}{rk}
\DeclareMathOperator{\mult}{mult}
\DeclareMathOperator{\id}{id}

\title[Cyclic sieving phenomena in the cluster complex]{Cyclic sieving phenomena on parabolic classes of faces of the cluster complex}
\author{Lucas Pouillart}
\address{Université Paris Cité, CNRS, IRIF, F-75013, Paris, France}
\date{}							

\begin{document}

\newtheorem{Thm}{Theorem}[section]
\newtheorem{Lm}[Thm]{Lemma}
\newtheorem{Prop}[Thm]{Proposition}
\newtheorem{Cor}[Thm]{Corollary}

\theoremstyle{definition}
\newtheorem{Def}{Definition}[section]
\newtheorem{ex}[Def]{Example}

\theoremstyle{remark}
\newtheorem*{rqe}{Remark}
\newtheorem*{pr}{Proof}
\newtheorem*{tm}{Toy model}

\begin{abstract}
The cyclic sieving phenomenon was introduced by Reiner, Stanton and White in 2004 as a generalization of Stembridge's $q=-1$ phenomenon. In a paper from 2008, Eu and Fu studied many occurrences of this phenomenon on the faces of the generalized cluster complex with the action of the Fomin-Reading rotation in the classical types $A_n$, $B_n$, $D_n$ and $I_2(k)$. There was yet no known uniform $q$-analogue of the $k$-face numbers of these complexes. In a more recent paper from 2023, Douvropoulos and Josuat-Vergès provided a refinement of the enumeration of the faces of the generalized cluster complex using a uniform formula. For a parabolic subgroup $W_X \subset W$ of the associated Coxeter group $W$, their formula factorises nicely under the assumption that $N_W(W_X)/W_X$ acts as a reflection group on $X$, which is very often the case. Using this condition, we provide a uniform refinement of these cyclic sieving phenomena using a $q$-analogue of their main formula with a type by type proof based on the classification of finite irreducible Coxeter groups.
\end{abstract}

\maketitle

\section{Introduction}

Let $X \circlearrowleft C$ be a finite set acted on by a finite cyclic group $C$. One says that together with $P \in \mathbb{Z}[q]$, the triple $(X,C,P(q))$ exhibits the \textit{cyclic sieving phenomenon} (CSP for short) if the following relation holds:
\[
P(e^{\frac{2i\pi}{d}}) = | X^g |
\]
where $g \in C$ has order $d$. In particular $P(1) = |X|$.\\
This phenomenon was introduced in \cite{RSW04} as a generalization of Stembridge's "$q = -1$" phenomenon and many of its occurences have been studied over the past two decades.\\
One of the main occurences of the CSP is the following theorem due to Reiner, Stanton and White in \cite{RSW04}:
\begin{Thm}\cite{RSW04}
Let $X$ be the set of triangulations of the $(n+2)$-gon using noncrossing diagonals. Let
\[
\Cat_n(q) := \displaystyle \frac{1}{[n+1]_q} \begin{bmatrix} 2n \\ n \end{bmatrix}_q.
\]
Let the group $C_{n+2}$ act on $X$ as the polygon rotation. Then $(X, C_{n+2}, \Cat_n(q))$ exhibits the CSP.
\end{Thm}
This theorem provides us with a CSP on a Catalan object, using the natural $q$-analogue of the classical Catalan numbers. The Catalan numbers have several generalizations, one being the \textit{Kirkman numbers} $\Kirk_{k,n}$, which count the dissection of the $(n+2)$-gon using $k$ noncrossing diagonals. Using a natural $q$-analogue of these numbers, Reiner, Stanton and White generalized their previous result:
\begin{Thm}\cite{RSW04}
Let $X$ be the set of dissections of a regular $(n+2)$-gon using $k$ noncrossing diagonals. Let
\[
\Kirk_{k,n}(q) := \displaystyle \frac{1}{[n+1]_q} \begin{bmatrix}
n+k+1 \\
k \end{bmatrix}_q
\begin{bmatrix}
n-1 \\
k\end{bmatrix}_q.
\]
Let the group $C_{n+2}$ act on $X$ as the polygon rotation. Then $(X, C_{n+2}, \Kirk_{k,n}(q))$ exhibits the CSP.
\end{Thm}
The dissection poset of the $(n+2)$-gon can be seen representation theoretically as the \textit{cluster complex} $\Gamma$ of type $A_{n-1}$, associated to the root system $A_{n-1}$. Its clockwise rotation is the \textit{Fomin-Zelevinsky rotation} of $\Gamma$ as defined in \cite{FZ03}. In \cite{EF08}, Eu and Fu exhibited several CSP on the faces of the cluster complex of irreducible types using natural $q$-analogues of counting formulas of the faces of a certain dimension. In particular they exhibited the following uniform CSP on the facets of irreducible cluster complexes using a natural $q$-analogue of the \textit{generalised Catalan numbers} $\Cat(\Phi)$:
\begin{Thm}\cite{EF08}
Let $\Phi$ be an irreducible root system with Coxeter number $h$ and exponents $e_1, \cdots, e_n$, let $X$ be the set of facets of $\Gamma(\Phi)$ acted on by the group $C_{h+2}$ as the Fomin-Zelevinsky rotation, and let
\[
\Cat(\Phi,q) := \prod_{i=1}^n \frac{[h + e_i +1]_q}{[e_i + 1]_q}.
\]
Then $(X, C_{h+2}, \Cat(\Phi,q))$ exhibits the cyclic sieving phenomenon.
\end{Thm}
One could hope to generalise the Kirkman numbers to other irreducible types as well, but the uniform formulas exhibited by Fomin and Reading in \cite{FR05} can only be expressed up to a "mysterious factor". Recent results by Douvropoulos and Josuat-Vergès in \cite{DJV23} provided us with the uniform counting formula $\mu_\lambda$ recalled in Theorem \ref{enumform} for faces of type $\lambda$. We can define a $q$-analogue $\mu_\lambda(q)$ of the formula $\mu_\lambda$ given by \eqref{rcandidate} in many cases. With this $q$-analogue, we get the main result of this paper:
\begin{Thm}
\label{main}
Let $\Gamma(\Phi)_\lambda$ be the set of faces of $\Gamma(\Phi)$ of parabolic type $\lambda$ and $\mathcal{R}$ be the Fomin-Zelevinsky rotation. Assume that there exists a parabolic subgroup $W_\lambda$ of type $\lambda$ such that $N(W_\lambda)/W_\lambda$ acts as a reflection group on $X := \Fix(W_\lambda)$ and let $\mu_\lambda(q)$ be the polynomial defined by \eqref{rcandidate}. Then $(\Gamma(\Phi)_\lambda, \langle \mathcal{R} \rangle, \mu_\lambda(q))$ exhibits the CSP.
\end{Thm}

The condition of $N_W(W_X)/W_X$ being a reflection group on $X$ is very common: it is automatic in types $A_n$, $B_n$, $I_2(k)$, $F_4$, $H_3$ and $H_4$, and is true in most cases in types $D_n$, $E_6$, $E_7$ and $E_8$. These groups have been studied thoroughly by Howlett \cite{H80} in the 80's, and a recent new description of their structure in \cite{DPR25} highlights the reason why this condition holds so often.

\begin{rqe}
Theorem \ref{main} partially solves an open problem given at the end of \cite{EF08}. CSPs on disjoint unions of $C$-sets behave nicely with the sieving polynomials: the sieving polynomials $P$ and $Q$ for $A$ and $B$ respectively give rise to a CSP on $A \uplus B$ with sieving polynomial $P + Q$. A natural $q$-analogue of the face numbers of the cluster complex is then given by an appropriate sum of polynomials $\mu_\lambda(q)$.
\end{rqe}

After some preliminaries about the generalized cluster complex in section \ref{ClusterComplex}, we will review the structure of normalizers of parabolic subgroups in Section \ref{normal} so that we can define our $q$-analogue in Section \ref{q-counting}. Sections \ref{PA}, \ref{PB}, \ref{PI} and \ref{PD} will be devoted to proving our main theorem for each classical family (namely $A_n$, $B_n$, $I_2(k)$ and $D_n$), and we will look at the case of exceptional root systems in Section \ref{exceptional}.

\section{The cluster complex}
\label{ClusterComplex}

Let $\Phi$ be a root system with reflection group $W$ and simple system $\Delta = \Delta_+ \sqcup \Delta_-$ such that $\Delta_\epsilon$ contains only roots which are pairwise orthogonal. Let $h$ be its  Coxeter number, and let $\Phi_{\geq -1} := -\Delta \cup \Phi_+$ be its set of \textit{almost positive roots}. Consider the \textit{bipartite Coxeter element} $ c := c_+ c_-$ where $c_\epsilon = \prod_{i \in \Delta_\epsilon} s_i$. Fomin and Zelevinsky defined in \cite{FZ03} the following action of $C_{h+2}$ on $\Phi_{\geq -1}$:
\begin{Def}\cite{FZ03}
The \textit{Fomin-Zelevinsky rotation} of $\Phi_{\geq -1}$ is the bijection $\mathcal{R} : \Phi_{\geq -1} \rightarrow \Phi_{\geq -1}$ defined by:
\[ \mathcal{R}(\alpha) =
\begin{cases}
- \alpha &\text{ if } \alpha \in (- \Delta_+) \cup \Delta_-,\\
c(\alpha) &\text{ otherwise.}
\end{cases}
\]
This rotation has order $\frac{h+2}{2}$ whenever $- \id \in W$ (namely in irreducible types $B_n$, $D_{2k}$, $I_2(2k)$, $E_7$, $E_8$, $F_4$, $H_3$ and $H_4$) and $h+2$ otherwise (in irreducible types $A_n$, $D_{2k+1}$, $I_2(2k+1)$ and $E_6$) and thus induces an action of $C_{h+2}$ on $\Phi_{\geq -1}^{(m)}$.
\end{Def}
This action can be generalized to the following object:
\begin{Def}\cite{FR05}
We call \textit{$m$-colored almost positive roots} and denote by $\Phi^{(m)}_{\geq -1}$ the union of the following subsets of $\Phi \times \mathbb{N}$:
\begin{itemize}
\item $\Phi^{(m)}_+ := \{\alpha^i \text{ : } \alpha \in \Phi_+, 1 \leq i \leq m\}$,
\item $\{ - \alpha^1 \text{ : } \alpha \in - \Delta\}$.
\end{itemize}
The \textit{Fomin-Reading rotation} $\mathcal{R}_m$ is the following bijection $\mathcal{R} : \Phi^{(m)}_{\geq -1} \rightarrow \Phi^{(m)}_{\geq -1}$:
\[ \mathcal{R}_m(\alpha^i) =
\begin{cases}
\alpha^{i+1} &\text{ if } \alpha^i \in \Phi_+^{(m)} \text{ and } i < m, \\
\mathcal{R}(\alpha^1) &\text{ otherwise.}
\end{cases}
\]
\end{Def}
In the case $m = 1$, we get the classical almost positive roots and the Fomin-Zelevinsky rotation. The Fomin-Reading rotation has order $d \mid mh + 2$ and induces an action of $C_{mh+2}$ on $\Phi^{(m)}_{\geq -1}$.\\
From this action we can define the following binary relation:
\begin{Def}\cite{FR05}
Let $||$ be the binary relation on $\Phi^{(m)}_{\geq -1}$ defined by the following conditions:
\begin{itemize}
\item if $\alpha \in \Delta$, then $- \alpha \text{ }|| \text{ } \beta$ if and only if $\alpha$ does not appear in the decomposition of $\beta$ as a linear combination of simple roots,
\item $\mathcal{R}_m(\alpha)\text{ }||\text{ } \mathcal{R}_m(\beta) \iff \alpha \text{ }||\text{ }\beta$.
\end{itemize}
This relation is called the \textit{compatibility relation} of $\Phi^{(m)}_{\geq -1}$.
\end{Def}
The \textit{$m$-cluster complex} $\Gamma^{(m)}$ of $\Phi$ is the flag complex of the compatibility relation on $\Phi^{(m)}_{\geq -1}$. In the case $m=1$ and $\Phi$ cristallographic, it has an interpretation in terms of the cluster variables of a cluster algebra. Another representation-theoretic interpretation in the general case is due to Thomas in \cite{T08} as a complex of irreducible components of cluster-tilting objects in a certain category.
\begin{tm}
The classical example of cluster complex is given by $\Gamma^{(m)}(A_{n-1})$, the $m$-cluster complex of type $A_{n-1}$. Type $A$ objects are often classical combinatorial objets and its cluster complexes make no exceptions. The type $A_{n-1}$ $m$-cluster complex can be realized as the complex of $m$-divisible dissections of the $(mn+2)$-gon using noncrossing diagonals. The action of the Fomin-Reading rotation on the faces of $\Gamma^{(m)}(A_{n-1})$ can be seen as the clockwise rotation of the $(mn+2)$-gon. In particular, its facets are the $(m+2)$-angulations of the $(mn+2)$-gon, and are counted by the famous \textit{Fuß-Catalan number} $\Cat_n^m$.
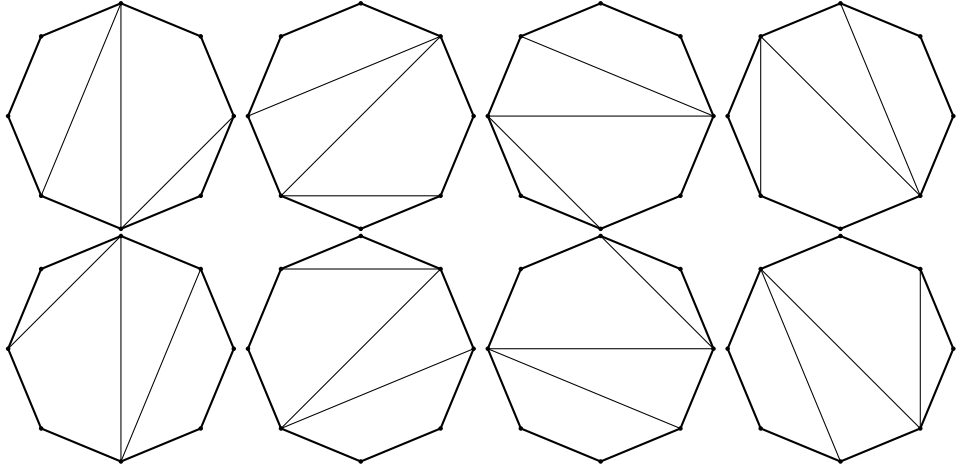
\begin{figure}[htbp]
\begin{tikzpicture}[scale=1.5]
  \def\n{8}
  \def\r{1}
  \foreach \i in {1,...,\n} {
    \coordinate (P\i) at ({\r*cos(90 - 360/\n*(\i-1))}, {\r*sin(90 - 360/\n*(\i-1))});
  }
  \foreach \i in {1,...,\n} {
    \pgfmathtruncatemacro{\j}{mod(\i,\n)+1}
    \draw[thick] (P\i) -- (P\j);
  }
  \draw (P1) -- (P5);
  \draw (P3) -- (P5);
  \draw (P1) -- (P6);
  \foreach \i in {1,...,\n} {
    \fill (P\i) circle (0.02);
  }
\end{tikzpicture}
\begin{tikzpicture}[scale=1.5]
  \def\n{8}
  \def\r{1} 
  \foreach \i in {1,...,\n} {
    \coordinate (P\i) at ({\r*cos(90 - 360/\n*(\i-1))}, {\r*sin(90 - 360/\n*(\i-1))});
  }
  \foreach \i in {1,...,\n} {
    \pgfmathtruncatemacro{\j}{mod(\i,\n)+1}
    \draw[thick] (P\i) -- (P\j);
  }
  \draw (P2) -- (P6);
  \draw (P4) -- (P6);
  \draw (P2) -- (P7);
  \foreach \i in {1,...,\n} {
    \fill (P\i) circle (0.02);
  }
\end{tikzpicture}
\begin{tikzpicture}[scale=1.5]
  \def\n{8}
  \def\r{1}
  \foreach \i in {1,...,\n} {
    \coordinate (P\i) at ({\r*cos(90 - 360/\n*(\i-1))}, {\r*sin(90 - 360/\n*(\i-1))});
  }
  \foreach \i in {1,...,\n} {
    \pgfmathtruncatemacro{\j}{mod(\i,\n)+1}
    \draw[thick] (P\i) -- (P\j);
  }
  \draw (P3) -- (P7);
  \draw (P5) -- (P7);
  \draw (P3) -- (P8);
  \foreach \i in {1,...,\n} {
    \fill (P\i) circle (0.02);
  }
\end{tikzpicture}
\begin{tikzpicture}[scale=1.5]
  \def\n{8}
  \def\r{1}
  \foreach \i in {1,...,\n} {
    \coordinate (P\i) at ({\r*cos(90 - 360/\n*(\i-1))}, {\r*sin(90 - 360/\n*(\i-1))});
  }
  \foreach \i in {1,...,\n} {
    \pgfmathtruncatemacro{\j}{mod(\i,\n)+1}
    \draw[thick] (P\i) -- (P\j);
  }
  \draw (P4) -- (P8);
  \draw (P6) -- (P8);
  \draw (P4) -- (P1);
  \foreach \i in {1,...,\n} {
    \fill (P\i) circle (0.02);
  }
\end{tikzpicture}
\begin{tikzpicture}[scale=1.5]
  \def\n{8}
  \def\r{1}
  \foreach \i in {1,...,\n} {
    \coordinate (P\i) at ({\r*cos(90 - 360/\n*(\i-1))}, {\r*sin(90 - 360/\n*(\i-1))});
  }
  \foreach \i in {1,...,\n} {
    \pgfmathtruncatemacro{\j}{mod(\i,\n)+1}
    \draw[thick] (P\i) -- (P\j);
  }
  \draw (P5) -- (P1);
  \draw (P7) -- (P1);
  \draw (P5) -- (P2);
  \foreach \i in {1,...,\n} {
    \fill (P\i) circle (0.02);
  }
\end{tikzpicture}
\begin{tikzpicture}[scale=1.5]
  \def\n{8}
  \def\r{1}
  \foreach \i in {1,...,\n} {
    \coordinate (P\i) at ({\r*cos(90 - 360/\n*(\i-1))}, {\r*sin(90 - 360/\n*(\i-1))});
  }
  \foreach \i in {1,...,\n} {
    \pgfmathtruncatemacro{\j}{mod(\i,\n)+1}
    \draw[thick] (P\i) -- (P\j);
  }
  \draw (P6) -- (P2);
  \draw (P8) -- (P2);
  \draw (P6) -- (P3);
  \foreach \i in {1,...,\n} {
    \fill (P\i) circle (0.02);
  }
\end{tikzpicture}
\begin{tikzpicture}[scale=1.5]
  \def\n{8}
  \def\r{1}
  \foreach \i in {1,...,\n} {
    \coordinate (P\i) at ({\r*cos(90 - 360/\n*(\i-1))}, {\r*sin(90 - 360/\n*(\i-1))});
  }
  \foreach \i in {1,...,\n} {
    \pgfmathtruncatemacro{\j}{mod(\i,\n)+1}
    \draw[thick] (P\i) -- (P\j);
  }
  \draw (P7) -- (P3);
  \draw (P1) -- (P3);
  \draw (P7) -- (P4);
  \foreach \i in {1,...,\n} {
    \fill (P\i) circle (0.02);
  }
\end{tikzpicture}
\begin{tikzpicture}[scale=1.5]
  \def\n{8}
  \def\r{1}
  \foreach \i in {1,...,\n} {
    \coordinate (P\i) at ({\r*cos(90 - 360/\n*(\i-1))}, {\r*sin(90 - 360/\n*(\i-1))});
  }
  \foreach \i in {1,...,\n} {
    \pgfmathtruncatemacro{\j}{mod(\i,\n)+1}
    \draw[thick] (P\i) -- (P\j);
  }
  \draw (P8) -- (P4);
  \draw (P2) -- (P4);
  \draw (P8) -- (P5);
  \foreach \i in {1,...,\n} {
    \fill (P\i) circle (0.02);
  }
\end{tikzpicture}

\caption{The orbit under the rotation of a face of $\Gamma(A_{5})$.}
\end{figure}
\end{tm}
Instead of classifying the faces of $\Gamma^{(m)}(W)$ by dimension, we will use a thinner criterion, its \textit{parabolic type}. Let us first recall some definitions and facts about generalized noncrossing partitions:

\begin{Def}
The \textit{intersection lattice} of $\Phi$, which we denote by $\Pi(\Phi)$ or $\Pi(W)$ is the poset defined by the following:
\begin{itemize}
\item its elements are the intersections of the reflecting hyperplanes $H_\alpha$, $\alpha \in \Phi_+$,
\item $X \leq Y \iff X \supseteq Y$ for all $X,Y \in \Pi(W)$.
\end{itemize}
This poset is also called the \textit{generalised partition poset} of the root system.\\
If $X$ is a flat, we denote by $\Pi(W)_X$ the intersection lattice obtained by replacing the hyperplanes $H_{\alpha}$ by $(H_{\alpha} \cap X)$.
\end{Def}

\begin{rqe} We can notice the following:
\begin{itemize}
\item The lattice $\Pi(W)$ is isomorphic to the lattice of parabolic subgroups. Let $X \in \Pi(W)$ and $W_X$ be the pointwise stabilizer of $X$ in $W$. The map $X \mapsto W_X$ is a lattice isomorphism.
\item In type $A$, these lattices are the classical set partition posets. The other classical types also admit good combinatorial descriptions using certain set partitions with some symmetry conditions.
\item When $X \in \Pi(W)$, $\Pi(W)_X$ is naturally isomorphic to the interval $[X, \{0\}]$ in $\Pi(W)$.
\end{itemize}
\end{rqe}

Let us denote by $T$ the set of all reflections of $W$, and $l_T(w)$ the \textit{reflection length} of $w \in W$. The \textit{absolute order} on $W$ is the partial order relation defined on $W$ by
\[
u \preceq v \iff l_T(v) = l_T(u) + l_T(u^{-1}v).
\]
The interval $[e,c]$ in the absolute order is called the \textit{generalised noncrossing partition lattice} of $W$ which we denote by $\NC(W)$. It is naturally isomorphic to a subposet of $\Pi(W)$ by $w \mapsto \Fix(w)$.

\begin{Prop}
For every face $f \in \Gamma^{(m)}(W)$, there exists an indexing $f = \{ \alpha_1^{i_1}, \cdots, \alpha_k^{i_k}\}$ such that
\[
\prod f := \prod_{i=1}^k t_{\alpha_i} \in \NC(W)
\]
where $t_{\alpha_i}$ is the reflection associated to the root $\alpha_i$.
\end{Prop}
This proposition is a consequence of an alternative definition of $\Gamma^{(m)}(W)$ due to Athanasiadis and Tzanaki in \cite{AT08}.\\
We can now define $\underline{f} := c_+(\prod f) c_-$ which is also a noncrossing partition of $W$.
\begin{rqe}
The application $w \mapsto c_+ w c_-$ is an antiautomorphism of $\NC(W)$, called the \textit{bipartite Kreweras complement}.
\end{rqe}

To each face $f \in \Gamma^{(m)}(W)$, we associate the noncrossing partition $\underline{f}$. As seen before, a noncrossing partition is a parabolic subgroup $W_{\underline{f}}$ in a certain way, so we can associate a parabolic subgroup to each face of the complex. The following theorems due to Douvropoulos and Josuat-Vergès in \cite{DJV23} justify this choice:

\begin{Thm} 
The link of $f \in \Gamma^{(m)}(W)$ is isomorphic to $\Gamma^{(m)}(W_{\underline{f}})$.
\end{Thm}

\begin{Thm}
For every $f \in \Gamma^{(m)}(W)$, $W_{\underline{f}}$ and $W_{\underline{\mathcal{R}_m(f)}}$ are parabolically conjugate.
\end{Thm}

In other terms, it is enough to study the action of the rotation in the parabolic classes of the cluster complex. We denote by $\Gamma^{(m)}(W)_\lambda$ the set of faces of $\Gamma^{(m)}(W)$ of parabolic type $\lambda$. The main result of \cite{DJV23} states the following:

\begin{Thm}
\label{enumform}
Assume that $W$ is irreducible, and let $W_X \subset W$ be a parabolic subgroup of type $\lambda$ the faces of $\Gamma^{(m)}(W)_\lambda$ are counted by:
\begin{equation}
\label{DJV}
\mu_\lambda := (-1)^{\dim X} \frac{p_X(-mh-1)}{[N(W_X):W_X]},
\end{equation}
where $p_X$ is the characteristic polynomial of $\Pi(W)_X$.
\end{Thm}
As these numbers refine the face numbers of the generalized cluster complex, which are the generalized Kirkman numbers, we call them the \textit{refined Kirman numbers}.

\section{Normalizers of parabolic subgroups}
\label{normal}

In this section, let $W$ be a finite Coxeter group, $V$ be the ambient vector space of its geometric representation and $P$ be a parabolic subgroup of $W$. The numbers $[ N_W(P): P ]$ appearing in \eqref{DJV} have appeared many times in the theory of reflection groups and hyperplane arrangements with a lot of other quotient groups involving the parabolic subgroups of $W$. Following \cite{H80} and \cite{DPR25}, let us study two decompositions of $N_W(P)$.

\subsection{Howlett's decomposition}
Before anything, let us introduce the notion of \textit{Howlett complement} of a group through the following lemma:
\begin{Lm} [\cite{H80}, Lemma 2]
Let $G$ be a group of orthogonal automorphisms of $V$ and $W$ be a reflection group in $V$ with root system $\Phi$ such that $W \mathrel{\unlhd} G$. Then $W$ has a complement $H$ (in the sense that $G = WH$), where $H = \{ g \in G \mid g(\Phi_+) = \Phi_+ \}$. $H$ is called the \textit{Howlett complement of $W$ in $G$}.
\end{Lm}
It is obvious that $P \mathrel{\unlhd} N_W(P)$, so this lemma has the following consequence:
\begin{Cor}
If $J \subset \Delta$, then $N(W_J)/W_J \cong \{ w \in W \mid w(J) = J \}$.
\end{Cor}
This result comes from the identification of $\{ w \in W \mid w(J) = J \}$ as the \textit{Howlett complement of $W_J$} by Howlett. The fact that $P$ has a complement in $N_W(P)$ is a straightforward consequence of this corollary. This decomposition gives us the following decomposition of $N_W(P)$:
\[
N_W(P) = P \rtimes W'.
\]
Consequently $W' \cong N_W(P)/P$. It appears quickly that $W'$ is 'almost' a reflection group: in \cite{H80}, Howlett exhibits a subgroup $W'' \subset W'$ which is generated by special elements of $W'$ called the $R$-elements see [\cite{H80}, Page 2] for more details.
\begin{Lm}
The $R$-elements are reflections in $\Fix(P)$ so $W''$ acts as a reflection group on $\Fix(P)$. Furthermore, $W'' \mathrel{\unlhd} W'$ and $W''$ admits a complement in $W'$.
\end{Lm}
With that decomposition we get
\begin{equation}
\label{Howlett}
N_W(P) = P \rtimes ( W'' \rtimes D)
\end{equation}
for some $D \subset W'$ which was explicitly described by Howlett [\cite{H80}, Corollary 7].\\
This decomposition already allowed Howlett to look at the action of $N_W(P)/P$ on $\Fix(P)$ in most cases, but let us look at another recent description of this group.
\subsection{Douglass-Pfeiffer-Röhrle decomposition}
This decomposition is really recent at the time we are writing this sentence. The idea behind this other decomposition is that $N_W(P)$ can actually be realized as a subgroup of $O(X) \times O(X^{\perp})$, where $X = \Fix P$ so $P$ acts on $X^\perp$. The interesting part for our work is that it allows us to decompose the action of $N_W(P)$ into a direct sum of three $N_W(P)$-modules, two of which are of direct interest for us. Let $T_W$ be the set of reflections of a reflection group $W$.
\begin{Def}
Let $U$ be a reflection subgroup of $W$, the \textit{orthogonal complement} of $U$ is the group $U^{\dagger}$ defined by
\[
U^{\dagger} = \langle t \in T_W \mid ts = st \text{ for all } s \in U \cap T_W \rangle.
\]
\end{Def}
This definition is crucial for our decomposition: let $Q = P^{\dagger}$. It is obvious that:
\begin{itemize}
\item $Q$ acts on $X$,
\item $P \times Q \subset N_W(P)$,
\item $Q$ is a parabolic subgroup of $W$, which allows us to decompose $V = Y \oplus Y^\perp$ with $Y = \Fix Q$.
\end{itemize}
As $X^{\perp} \cap Y^{\perp} = \{ 0 \}$ and $Y^{\perp} \subset X$, we get that $V = X^{\perp} \oplus (X \cap Y) \oplus Y^\perp$.

\begin{Prop}
$P \times Q$ has a Howlett complement $D$ in $N_W(P)$.
\end{Prop}

The idea of the rest of the decomposition is to use Goursat's lemma [\cite{DPR25}, Section 4] to identify $D$ as a subgroup of $O(X^\perp) \times O(X \cap Y) \times O(Y^{\perp})$. It is naturally described in Sections 5 and 6 of \cite{DPR25} as a product $D = (A \times B) \rtimes C$ using Goursat's isomorphism [\cite{DPR25}, Section 4]. $(A \times B)$ can always be complemented because $C$ is often trivial, and we can explicitly find $D = C$ case by case otherwise.

\begin{Prop}
$X = \Fix(P) = (X \cap Y) \oplus Y^{\perp}$ and the decomposition $N_W(P) = (P \times Q) \rtimes ((A \times B) \rtimes C)$ is such that:\begin{itemize}
\item $P$ fixes $X$,
\item $Q$ fixes $Y = (X \cap Y) \oplus X^{\perp}$,
\item $A$ fixes $Y^{\perp}$,
\item $B$ fixes $X \cap Y$.
\end{itemize}
In other words, we get that:
\begin{itemize}
\item $X^{\perp}$ is a $(P \rtimes ((A \times B) \rtimes C))$-module,
\item $X \cap Y$ is a $(A \rtimes C)$-module,
\item $Y^\perp$ is a $(Q \rtimes (B \rtimes C))$-module.
\end{itemize}
\end{Prop}
It is too much to ask $N_W(P)$ to be a reflection group on $V$, but this decomposition helps us look clearly at the $W'$ of Howlett's decomposition to see if $N_W(P)/P$ is a reflection group. As $P$ only acts nontrivially on $X^\perp$, we need to look at the action of $N_W(P)$ on $X = (X \cap Y) \oplus Y^\perp$. This is done precisely for every irreducible group in \cite{DPR25}.

\begin{Cor}
The group $N_W(P)/P$ is a reflection group on $X$ if and only if $C$ is trivial, $A$ is a reflection group on $X \cap Y$, and $Q \rtimes B$ is a reflection group on $Y^{\perp}$.
\end{Cor}

In what follows, we will give combinatorial descriptions of $N_W(P)/P$ for $P \subset W$, where $W$ is a group of a classical type.

\subsection{Parabolic quotients in type $A_{n-1}$}
\label{Adenom}

Let us recall some well-known facts: $A_{n-1} = \mathfrak{S}_n$ and its Coxeter number is $n$. The parabolic subgroups of $\mathfrak{S}_n$ are in bijection with the partitions of the set $[n]$:
\begin{itemize}
\item consider a partition $P = \{ B_1, \cdots, B_k \}$ of $[n]$, the parabolic subgroup of $\mathfrak{S}_n$ associated with $P$ is $\mathfrak{S}_{B_1} \times \cdots \times \mathfrak{S}_{B_k}$,
\item two partitions $P_1$ and $P_2$ give rise to two parabolically conjugate subgroups if and only if their block sizes form the same integer partition.
\end{itemize}
As a consequence, the parabolic conjugacy relation on Young subgroups is exceptionally exactly the isomorphism relation.\\
When we consider the representation of $\mathfrak{S}_n$ as permutation matrices, a parabolic subgroup $W_X = \prod \mathfrak{S}_k^{i_k}$ is isomorphic to a subgroup of block-diagonal matrices. These matrices contain one $k$-size permutation block for every factor of type $\mathfrak{S}_k$ of $W_X$, or equivalently a $k$-size permutation block for every part $k$ of the associated $\lambda \vdash n$. In that case, $N(W_X)/W_X$ is isomorphic to the group of permutations of same-size blocks. In other terms, for a parabolic subgroup $W_X = \mathfrak{S}_{1}^{n_1} \times \cdots \times \mathfrak{S}_{k}^{n_k}$ (so $n_i = \mult(i,\lambda)$) of rank $n-1-l$, we have:
\[
N(W_X)/W_X \cong A_{n_1-1} \times \cdots \times A_{n_k-1}.
\]
This is explicitly proved in \cite{H80} and \cite{DPR25}, along with the fact that this group acts as a reflection group on $X$.

\subsection{Parabolic quotients in type $B_n$}
\label{Bdenom}

Let us recall that $B_n$ is the $n$-th group of signed permutations, its Coxeter number is $2n$ and its parabolic subgroups may have type $A$ and $B$ irreducible components. The parabolic conjugacy relation in type $B_n$ has been described by Geck and Pfeiffer in \cite{GP00} in the same way as in type $A$:
\begin{itemize}
\item If $W_X$ does not contain any transposition of the form $(j \text{ }-j)$, then $W_X \subset A_{n-1} \subset B_n$ form a chain of parabolic subgroups. We can associate the same $\lambda \vdash n$ as in type $A_{n-1}$.
\item If it does contain such a transposition, then $W_X$ has a type $B_m$ component. One can associate a partition $\lambda \vdash n-m$ describing the type $A$ components as a parabolic subgroup of $A_{n-m-1}$ using the chain $W_X / B_m \subset A_{n-m-1} \subset B_n$.
\end{itemize}

\begin{rqe}
It makes sense combinatorially to consider the diagram $B_1 \cong A_1$ as the only subdiagram of $B_n$ containing the remarkable vertex.
\end{rqe}

The parabolic quotients $N(W_X)/W_X$ behave similarly as in type $A$. They act as reflection groups on $X$, see \cite{H80}. We give an explicit representation of $N_W(W_X)$ for every $W_X \subseteq B_n$ which allows us to describe explicitly its quotient by $W_X$. We can use the representation of $B_n$ as the group of $n \times n$ signed permutation matrices. This representation gives us a representation of a parabolic subgroup $W_X$ of type $\lambda \vdash k \leq n$ which can be described as the following group of block-diagonal matrices:
\begin{itemize}
\item For each block size $k$ of $\lambda$, there is a $k \times k$ sized block corresponding to classical permutations matrices.
\item If the type $B_j$ component exists, then we have a $j \times j$ sized block of signed permutation matrices.
\end{itemize}
These representations are faithful, so $W_X$ is isomorphic to this group. Up to parabolic conjugacy, $W_X$ can be illustrated by the following:
\[
\begin{bmatrix}
    A_{i_1}       &  &  &  \\
           & \ddots & & \\
     & & A_{i_k} & \\
      &  &  & B_j
\end{bmatrix}.
\]
Using this representation, let us look at $N_W(W_X)$:
\begin{itemize}
\item The type $A_k$ blocks are preserved only by signed permutation matrices exchanging type $A_k$ blocks with coefficients of the same sign on each block. 
\item If the type $B_j$ component exists, it is preserved by any action of a type $B_j$ matrix. 
\end{itemize}
At the quotient level, the type $B_j$ block vanishes and the only remaining information remaining on the type $A_k$ blocks are the blocks permuted and the sign of the permutation.
For $W_X$ of type $\lambda$ where for an integer $i$, $n_i = \mult(i, \lambda)$, we get the following:
\[
N(W_X)/W_X \cong B_{n_1} \times \cdots \times B_{n_k}.
\]

\subsection{Parabolic quotients in type $D_n$}
\label{Ddenom}

Let us recall that $D_n$ is the index 2 subgroup of $B_n$ defined by
\[
\sigma \in D_n \iff | \{ i \in [n] \mid \sigma(i) < 0 \} | \text{ is even.}
\]
Its Coxeter number is $2n-2$, and its parabolic subgroups are index 1 or 2 subgroups of parabolic subgroups of $B_n$, containing type $A$ irreducible components in the same way, and at most one type $D$ component. The parabolic conjugacy classes has again been described by Geck and Pfeiffer in \cite{GP00} as integer partitions in the following way:
\begin{itemize}
\item if $W_X$ has a type $D_k$ component, then $W_X$ gives rise to the same $\lambda \vdash n-k$ as in type $B_n$.
\item if $W_X$ has no type $D$ component, then we associate the same $\lambda \vdash n$ as in type $B_n$ but the classification is not over:
\begin{itemize}
\item if $\lambda$ has an odd part, the parabolic conjugacy class is entirely given by $\lambda$.
\item If $\lambda$ is all-even, then up to conjugacy exactly one of the two remarkable reflections $s_0$ or $s_1$ is in $W_X$. Thus we associate $(\lambda, \pm)$ to $W_X$ depending on which one belongs in $W_X$ up to conjugacy.
\end{itemize}
\end{itemize}
\begin{rqe}
In a combinatorial point of view, it makes sense to consider the diagrams $D_2 \cong A_1 \times A_1$ and $D_3 \cong A_3$ as subdiagrams of $D_n$ containing the remarkable vertices 0 and 1.
\end{rqe}
Now let $W_X \subset D_n$ be a parabolic subgroup of type $\lambda \vdash k \leq n$. If $k \leq n$, let $n_i = \mult(i,\lambda)$. Then if $k \leq n$, the $D$ factor in Howlett's decomposition \eqref{Howlett} is trivial, therefore $N(W_X)/W_X$ acts as a reflection group and we get:
\[
N(W_X)/W_X = \prod_{i \in \lambda} B_{n_i}.
\]
If $k = n$, then Howlett's decomposition is a bit more complicated. Let $I$ be the number of odd sizes of parts in $\lambda$, then
\[
W'' = \prod_{i \text{ odd}} D_{n_i} \times \prod_{i \text{ even}} B_{n_i} \text{ and } D = \mathbb{Z}_2^{I - 1}.
\]
It is immediate to see that if $I = 1$, then $W' = W''$ and it is therefore a reflection group on $X$. If $I > 1$, then $W'$ is not a reflection group on $X$. This could be deduced from \cite{H80} and \cite{DPR25} but we will give a way simpler argument.
Parabolic subgroups of $D_n$ can, just like in type $B_n$, be represented as the following group of block-diagonal $n \times n$-matrices:
\begin{itemize}
\item For each block size $k$ of $\lambda$, there is a $k \times k$ sized block corresponding to classical permutations matrices.
\item If the type $D_j$ component exists, then we have a $j \times j$ sized block of signed permutation matrices.
\end{itemize}
These representations are also faithful, so $W_X$ is isomorphic to this group. Up to parabolic conjugacy, $W_X$ can be illustrated by the following:
\[
\begin{bmatrix}
    A_{i_1}       &  &  &  \\
           & \ddots & & \\
     & & A_{i_k} & \\
      &  &  & D_j
\end{bmatrix}.
\]
In that case $N(W_X)/W_X$ is a subgroup of index 1 or 2 of the corresponding group in type $B_n$. Two even block sizes can be permuted at will, but only an even number of pairs of odd block sizes can be permuted in one go. Then it is an index 2 subgroup if and only if $\lambda$ has an odd part. The only possible way to obtain an index 2 subgroup of $B_{n_1} \times \cdots \times B_{n_k}$ with the structure of a reflection group is to reduce one of the $B_l$ component to $D_l$. it is obvious that two such groups are not isomorphic if $I \neq 1$.

\subsection{Parabolic quotients in dihedral type}

Let us consider here the dihedral type $W = I_2(k) = \langle s_1, s_2 \rangle$. This family contains exactly the rank 2 Coxeter groups. There are exactly 4 standard parabolic subgroups in these cases, two of which are the trivial parabolic subgroups ($\{ e \}$ and $I_2(k)$). The only two remaining cases, namely $\langle s_1 \rangle$ and $\langle s_2 \rangle$ are isomorphic of type $A_1$. The problem of deciding whether they are parabolically conjugate or not is well known: all rank 1 parabolic subgroups of $I_2(k)$ are conjugate if and only if $k$ is odd. Even if they are not conjugate, their normalizers are isomorphic because of the automorphism defined by $s_1 \leftrightarrow s_2$. The case of these groups is explicitly treated in [\cite{DPR25}, Tables 12 and 13).

\begin{Prop}
 Let $k \geq 3$ and $W_X = \langle t_\alpha \rangle$:
\begin{itemize}
\item if $k$ is odd, $N_W(W_X)/ W_X \cong A_0$,
\item if $k$ is even, $N_W(W_X)/W_X \cong A_1$.
\end{itemize}
\end{Prop}

\begin{rqe}
In this case, $N(W_X)/W_X$ is always a reflection group on $X$: in rank 0 and 2 this is obvious, and in rank 1 $A_0$ is a reflection group. It is enough to notice that $\dim X = 1$ and $A_1 \cong \mathbb{Z}_2$ can only act faithfully as the reflection group generated by $1 \leftrightarrow -1$.
\end{rqe}

\section{$q$-analogues of counting formulas}
\label{q-counting}

From now on, let $W$ be an irreducible Coxeter group of Coxeter number $h$. The previous formula \eqref{DJV}
\[
\mu_X := (-1)^{\dim X} \frac{p_X(-mh-1)}{[N(W_X):W_X]}
\]
has a natural $q$-analogue in many cases, let us describe this $q$-analogue in this section.\\
Let us first notice that the polynomial $p_X$ splits, and that its roots, called the \textit{exponents} of the lattice $\Pi(W)_X$, are positive integers. The exponents of every $\Pi(W)_X$, with $W$ finite irreducible, have all been computed by Orlik and Solomon in \cite{OS83}. The $(-1)^{\dim X}$ factor is here to ensure that the result is positive, so we can rewrite the numerator of \eqref{DJV} in the following way:
\begin{equation}
\label{num}
(-1)^{\dim X} p_X(-mh-1) = \prod_{i=1}^k (e_i^X + 1 + mh),
\end{equation}
where $e_i^X$ are the exponents of $\Pi(W)_X$. As mentioned above, they have been explicitly computed by Orlik and Solomon in many cases in \cite{OS83}. Let us look at some examples:
\begin{itemize}
\item The type $A_n$ partition lattice is isomorphic to the classical partition lattice of the set $[n+1]$ .
\item The type $B_n$ partition lattice is isomorphic to the $n$-th signed partition lattice.
\item The type $D_n$ partition lattice is isomorphic to the sublattice of $B_n$ obtained by forbidding the blocks of shape $\{i,-i\}$.
\item The type $I_2(k)$ partition lattice is the rank 2 lattice with unique minimal and maximal elements and $k$ rank 1 elements, see Figure 2.
\end{itemize}
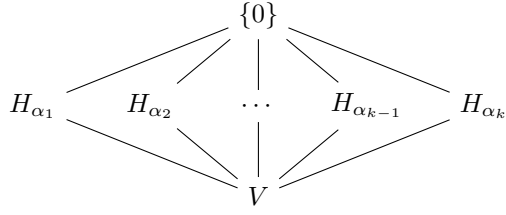
\begin{figure}[h]
\label{I-lattice}
\begin{tikzpicture}[scale=1.2, every node/.style={font=\small}]
  \node (zero) at (0,2) {$\{0\}$};
  \node (R2) at (0,0) {$V$};
  \node (H1)   at (-2.5,1) {$H_{\alpha_1}$};
  \node (H2)   at (-1.2,1) {$H_{\alpha_2}$};
  \node (dots) at (0,1)    {$\cdots$};
  \node (Hk1)  at (1.2,1)  {$H_{\alpha_{k-1}}$};
  \node (Hk)   at (2.5,1)  {$H_{\alpha_k}$};
  \foreach \x in {H1,H2,dots,Hk1,Hk}{
    \draw (zero) -- (\x);
    \draw (\x) -- (R2);
  }

\end{tikzpicture}
\caption{The intersection lattice of type $I_2(k)$}
\end{figure}
We need to look at the upper sets of $\Pi(W)$. This problem is easily solved in types $A_n$ and $B_n$, as they exhibit the same heredity property.
\begin{Prop}
\label{heredity}
If $W$ is an irreducible Coxeter group of type $A_n$ (\textit{resp.} $B_n$) and $X \in \Pi(W)$, the upper set $[X, \{0\}]$ is isomorphic to the lattice $\Pi(A_{n-\rk(X)})$ (\textit{resp.} $\Pi(B_{n-\rk(X)})$).
\end{Prop}
It is obvious that a strict upper set of $I_2(k)$ has either type $A_1$ or $A_0$ depending on its rank, so it remains to look at the case of type $D_n$. Say that an upper set of $\Pi(D_n)$ has type $D_n^k$ if it is generated by a partition $X$ with no central block and $k$ blocks of size $>1$. The exponents of such a lattice have been again computed in \cite{OS83}. It is useful to notice that $D_l^0 \cong D_l$ and $D_l^l \cong B_l$.
\begin{Prop} Let $W_X \subset D_n$ be a parabolic subgroup such that $\rk(W_X) = n-l$.
\begin{itemize}
\item If $W_X$ has type $\lambda \vdash k < n$, then $\Pi(D_n)_X$ has type $B_l$.
\item If $W_X$ has type $\lambda \vdash n$ with no 1-part in $\lambda$, then $\Pi(D_n)_X$ has type $B_l$.
\item If $W_X$ has type $\lambda \vdash n$ with $i$ parts $\neq 1$ in $\lambda$, then $\Pi(D_n)_X$ has type $D_l^i$.
\end{itemize}
\end{Prop}

As for the denominator, let us first recall that if $G$ is a reflection group (on any field) then $|G|$ can be factorized as the product of its \textit{degrees} as defined with the Chevalley-Shephard-Todd theorem. The group $N(W_X)/W_X$ acts faithfully on $X$ as $X = \Fix W_X$, meaning that $W_X = \ker \rho_{| N(W_X)}$, where $\rho$ is the geometric representation of $W$. The decompositions of \cite{H80} and \cite{DPR25} reviewed in Section \ref{normal} show that in a lot of cases, the group $N(W_X)/W_X$ acts as a reflection group on $X$. In that case, as mentioned before we have:
\begin{equation}
\label{denom}
[N(W_X):W_X] = \prod_{i=1}^l d_i^X,
\end{equation}
where $d_i^X$ are the \textit{degrees} of $N(W_X)/W_X$.
By combining \eqref{num} and \eqref{denom}, we get:
\begin{equation}
\label{rcandidate}
\mu_\lambda = \frac{\prod_{i=1}^k (e_i^X + 1 + mh)}{\prod_{i=1}^l d_i^X}.
\end{equation}
A natural candidate for a $q$-analogue of this formula is the following:
\begin{equation}
\label{qana}
\mu_\lambda(q) := \frac{\prod_{i=1}^k [e_i^X + 1 + mh]_q}{\prod_{i=1}^l [d_i^X]_q}.
\end{equation}

\begin{rqe}
Whenever $N(W_X)/W_X$ acts on $X$ as a reflection group, the faithfulness of this action automatically gives us $l \leq \dim X$. As $k = \dim X$, we notice that in most cases we get $k = l$. The only cases where we actually get $k > l$ are found in type $E_6$ and $H_4$.
\end{rqe}

Throughout the proof of our main theorem, we will repeatedly use the following results:

\begin{Prop}
\label{q-prop}
Let $n$, $k$ and $d$ be positive integers, and let $\zeta_d$ be a primitive $d$-th root of unity. Then
\begin{itemize}
\item $[n]_{q = \zeta_d} = 0 \iff d \mid n$,
\item if $n \equiv k \mod d$, $\displaystyle \lim_{q \rightarrow \zeta_d} \frac{[n]_q}{[k]_q} = \begin{cases} \frac{n}{k} \text{ if } n \equiv k \equiv 0 \mod d \\ 1 \text{ otherwise,} \end{cases}$
\item ($q$-Lucas theorem) if $n = ad + b$ and $k = rd+s$, where $ 0 \leq b,s < d$, then
\[
\begin{bmatrix} m \\ k \end{bmatrix}_{q = \zeta_d} = \binom{a}{r} \begin{bmatrix} b \\ s \end{bmatrix}_{q = \zeta_d}.
\]
\end{itemize}
\end{Prop}

Now that $\mu_\lambda(q)$ is defined, we have all the tools to start proving our main theorem.

\section{Proof of the main theorem in type $A_{n-1}$}
\label{PA}

 A face of $\Gamma^{(m)}(A_{n-1})$ of type $\prod \mathfrak{S}_k^{i_k}$ is a dissection of the $(mn+2)$-gon containing exactly $i_k$ inner $(mk+2)$-gons for every $k$. One can associate a partition $\lambda \vdash n$ to any dissection: if a $(n+2)$-gon can be dissected into a family $S$ of subgons, then $\sum_{\tau \in S} k(\tau) = n$, where $k(\tau)+2$ is the number of edges of $\tau$. The set of $k(\tau)$'s is an integer partition who matches the one associated to the parabolic type of the dissection.

\begin{figure}[h]
\begin{tikzpicture}[scale=2]
  \def\n{8}
  \def\r{1}
  \foreach \i in {1,...,\n} {
    \coordinate (P\i) at ({\r*cos(90 - 360/\n*(\i-1))}, {\r*sin(90 - 360/\n*(\i-1))});
  }
  \foreach \i in {1,...,\n} {
    \pgfmathtruncatemacro{\j}{mod(\i,\n)+1}
    \draw[thick] (P\i) -- (P\j);
  }
  \draw (P4) -- (P8);
  \draw (P5) -- (P8);
  \draw (P8) -- (P2);
  \foreach \i in {1,...,\n} {
    \fill (P\i) circle (0.02);
  }
\end{tikzpicture}
\caption{A face of parabolic type $A_1^2$ (or $[1,1,2,2]$) of $\Gamma(A_5)$.}
\end{figure}
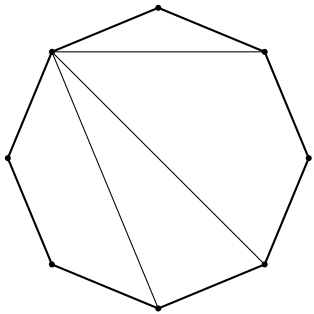

We will begin our proof by identifying our polynomial in type $A_{n-1}$. By Proposition~\ref{heredity}, the upper sets of the lattice can be identified as partition lattices themselves.

\begin{Prop}\cite{OS83}
The exponents of the lattice $\Pi(A_{n-1})_X$ are $1,2, \cdots, \dim X$.
\end{Prop}

Combining Subsection \ref{Adenom} and \cite{H80} gives us the set of $d_i^X$ for every $X \in \Pi(A_{n-1})$. For simplicity let us denote by $l$ the sum $n_1 + n_2 + \cdots + n_k$. We can now compute \eqref{num} for every parabolic subgroup $W_X$ of $A_{n-1}$. Following the notations of Section \ref{qana}, we get the polynomial:
\begin{align*}
\mu_\lambda(q) &=  \displaystyle \frac{\prod_{i=1}^l [e_i^X + 1 + n]_q}{\prod_{i=1}^l [d_i^X]_q},\\
&= \frac{\prod_{i=2}^k [n+i]_q}{[n_1]_q! [n_2]_q! \cdots [n_k]_q!}, \\
&= \frac{1}{[n+1]_q} \frac{\prod_{i=1}^k [n+i]_q}{[n_1]_q! [n_2]_q! \cdots [n_k]_q!}, \\
&= \frac{1}{[n+1]_q} \frac{[n+l]_q!}{[n_1]_q! [n_2]_q! \cdots [n_k]_q! [n]_q!}.
\end{align*}
The second factor can be rewritten as a $q$-multinomial coefficient, so we get:
\begin{equation}
\label{typeA}
\mu_\lambda(q) = \displaystyle \frac{1}{[n+1]_q} \begin{bmatrix} n+l \\ n_1, \cdots, n_k, n \end{bmatrix}_q.
\end{equation}
\begin{Thm}
The triple $(\Gamma(A_{n-1})_X, C_{n+2}, \mu_\lambda(q))$ exhibits the cyclic sieving phenomenon.
\end{Thm}

\begin{pr}
This is equivalent to Theorem 1 in \cite{AB25}:\begin{itemize}
\item we are counting dissections with $n_k$ inner $(k + 2)$-gons for each $k$,
\item these dissections obviously use $l-1$ noncrossing diagonals,
\item our polynomial $\mu_\lambda$ can be rewritten again as 
\[
\mu_\lambda(q) = \displaystyle \frac{1}{[n+1]_q}  \begin{bmatrix} n+l \\ l \end{bmatrix}_q \begin{bmatrix} l \\ n_1, \cdots, n_k \end{bmatrix}_q
\]
which gives us exactly the polynomial considered by Adams and Banaian for our case.
\end{itemize}
\end{pr}

\begin{rqe}
This cyclic sieving phenomenon also holds for the generalized cluster complexes of type $A$. The $C_{mn+2}$-set $\Gamma^{(m)}(A_{n-1})_X$ is isomorphic to $\Gamma(A_{mn-1})_Y$ where $W_Y \subset A_{mn-1}$ has an irreducible component of type $A_{km}$ for every irreducible component of type $A_k$ of $W_X$.
\end{rqe}

\section{Proof of the main theorem in type $B_n$}
\label{PB}

\subsection{Combinatorics of the type $B_n$ cluster complexes}

We can now start looking at the case of type $B_n$. The complex $\Gamma^{(m)}(B_n)$ can be combinatorially interpreted as the complex of $m$-divisible centrally symmetric dissections of the $(2mn+2)$-gon in the following way:
\begin{itemize}
\item A \textit{type $B_n$ diagonal} is either a diameter or a centrally symmetric pair of diagonals.
\item Two type $B_n$ diagonals are compatible if they are noncrossing and form only $(mk+2)$-gons for some integer $k$.
\item The Fomin-Reading rotation acts again as the classical clockwise rotation of the $(2mn+2)$-gon.
\end{itemize}
The parabolic type of a face $f \in \Gamma^{(m)}(B_n)$ can be found in the following way:\begin{itemize}
\item Consider $\mu = [1^{n_1}, 2^{n_2}, \cdots, k^{n_k}] \vdash 2n$ the partition associated to $f$ as a classical dissection.
\item Let $\lambda \vdash k \leq n$ be the integer partition $\lambda = [1^{\lfloor \frac{n_1}{2} \rfloor}, 2^{\lfloor \frac{n_2}{2} \rfloor}, \cdots, k^{\lfloor \frac{n_k}{2} \rfloor}]$.
\end{itemize}
The partition $\lambda$ characterizes the parabolic type of $f$ in the same fashion as the parabolic subgroups (see Subsection \ref{Bdenom}).One could notice that a face of $\Gamma^{(m)}(B_n)$ has a type $B_k$ component in its parabolic type if and only if it has a central $(2mk+2)$-gon. The components of type $A_l$ count the different subgons on a half of the polygon.
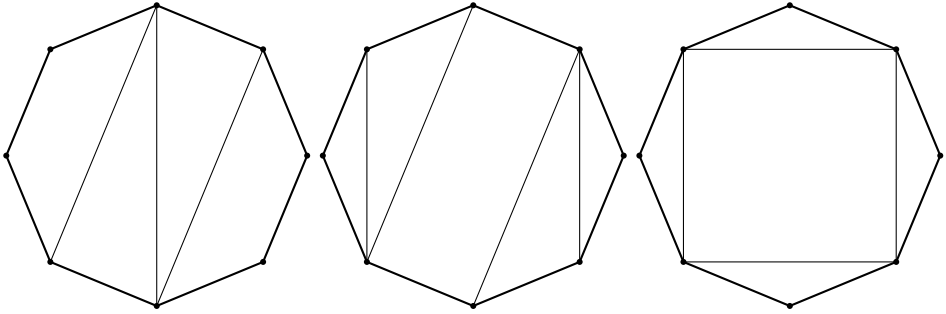
\begin{figure}[h]
\begin{center}
\begin{tikzpicture}[scale=2]
  \def\n{8}
  \def\r{1}
  \foreach \i in {1,...,\n} {
    \coordinate (P\i) at ({\r*cos(90 - 360/\n*(\i-1))}, {\r*sin(90 - 360/\n*(\i-1))});
  }
  \foreach \i in {1,...,\n} {
    \pgfmathtruncatemacro{\j}{mod(\i,\n)+1}
    \draw[thick] (P\i) -- (P\j);
  }
  \draw (P1) -- (P5);
  \draw (P2) -- (P5);
  \draw (P1) -- (P6);
  \foreach \i in {1,...,\n} {
    \fill (P\i) circle (0.02);
  }
\end{tikzpicture}
\begin{tikzpicture}[scale=2]
  \def\n{8}
  \def\r{1}
  \foreach \i in {1,...,\n} {
    \coordinate (P\i) at ({\r*cos(90 - 360/\n*(\i-1))}, {\r*sin(90 - 360/\n*(\i-1))});
  }
  \foreach \i in {1,...,\n} {
    \pgfmathtruncatemacro{\j}{mod(\i,\n)+1}
    \draw[thick] (P\i) -- (P\j);
  }
  \draw (P6) -- (P8);
  \draw (P2) -- (P4);
  \draw (P6) -- (P1);
  \draw (P5) -- (P2);
  \foreach \i in {1,...,\n} {
    \fill (P\i) circle (0.02);
  }
\end{tikzpicture}
\begin{tikzpicture}[scale=2]
  \def\n{8}
  \def\r{1}
  \foreach \i in {1,...,\n} {
    \coordinate (P\i) at ({\r*cos(90 - 360/\n*(\i-1))}, {\r*sin(90 - 360/\n*(\i-1))});
  }
  \foreach \i in {1,...,\n} {
    \pgfmathtruncatemacro{\j}{mod(\i,\n)+1}
    \draw[thick] (P\i) -- (P\j);
  }
  \draw (P8) -- (P2);
  \draw (P2) -- (P4);
  \draw (P4) -- (P6);
  \draw (P6) -- (P8);
  \foreach \i in {1,...,\n} {
    \fill (P\i) circle (0.02);
  }
\end{tikzpicture}
\end{center}
\caption{Faces of type $[1,2]$ on the left, and $[1,1]$ on the middle and the right, of $\Gamma(B_3)$.}
\end{figure}

\begin{rqe}
If $\lambda = [1^{n_1},2^{n_2}, \cdots, k^{n_k}]$ then $\Gamma^{(m)}(B_n)_\lambda \cong \Gamma^{(1)}(B_{mn})_{m\lambda}$ as $C_{2mn+2}$-sets.
\end{rqe}

\subsection{Sieving polynomials in type $B_n$}
\label{polyB}
Using the preceding remark, we treat the case $m=1$ and deduce the general case from it.
All the previous data allows us to compute the polynomial $\mu_\lambda(q)$ for a parabolic subgroup $W_X$ of type $\lambda = [1^{n_1},2^{n_2}, \cdots, k^{n_k}] \vdash t \leq n$ :
\begin{align*}
\mu_\lambda(q) &=  \displaystyle \frac{\prod_{i=1}^l [e_i^X + 1 + 2n]_q}{\prod_{i=1}^l [d_i]_q}\\
&= \frac{\prod_{i=1}^l [2n+2i]_q}{[2n_1]_q!! \cdots [2n_k]_q!!}\\
&= \frac{\prod_{i=1}^l [2n+2i]_q}{[2n_1]_q!! \cdots [2n_k]_q!!}.
\end{align*}
As all integers involved are even, this polynomial can be rewritten as:
\begin{align*}
\mu_\lambda(q) &= \frac{\prod_{i=1}^l [n+i]_{q^2}}{[n_1]_{q^2}! \cdots [n_k]_{q^2}!}\\
&= \begin{bmatrix} n+l \\ n_1, \cdots, n_k, n \end{bmatrix}_{q^2}.
\end{align*}
This expression of $\mu_X(q)$ is the one we are going to use to evaluate it on $\zeta_d$, where $\zeta_m$ is any $m$-th primitive root of unity.
\begin{Lm}
\label{evenB1}
If $d$ is odd, then $\mu_\lambda(\zeta_d) = \mu_\lambda(\zeta_{2d})$.
\end{Lm}
\begin{pr}
The proof here is very straightforward so we leave it as an exercise to the reader.
\end{pr}
\begin{rqe}
If $d$ is even and $d \mid 2n+2$, it is useful to notice that $\frac{d}{2} \mid n+1$, or equivalently $n \equiv (\frac{d}{2} -1) \mod \frac{d}{2}$.
\end{rqe}
A consequence of this lemma is that we only have to evaluate $\mu_X$ on even divisors of $2n+2$.
\begin{Lm}
\label{zeroB}
If $\mu_\lambda(\zeta_d) \neq 0$, then:
\begin{enumerate}
\item $\frac{d}{2} \mid l$,
\item $\frac{d}{2} \mid n_i$ for every $n_i$,
\item $d \mid 2j+2$, where $j = n-t$.
\end{enumerate}
\end{Lm}
\begin{pr}Let us check these case by case:
\begin{enumerate}
\item If $\frac{d}{2} \nmid l$, let us rewrite $\mu_\lambda(q)$ as the following expression:
\begin{equation}
\label{ABtypeB}
\mu_\lambda(q) = \begin{bmatrix} l \\ n_1, \cdots, n_k \end{bmatrix}_{q^2} \begin{bmatrix} n+l \\ n \end{bmatrix}_{q^2}.
\end{equation}
By $q$-Lucas theorem (Proposition \ref{q-prop}), when evaluated in $\zeta_d$, the second factor of the right-hand side of \eqref{ABtypeB} has a factor
\begin{equation}
\label{binom}
\begin{bmatrix} i \\ \frac{d}{2} - 1 \end{bmatrix}_{\zeta_d^2},
\end{equation}
with $i < \frac{d}{2}$. If $\frac{d}{2} \nmid l$, we get $i < \frac{d}{2} - 1$, which makes \eqref{binom} equal to 0.
\item In addition to what we previously assumed, let us assume that $\frac{d}{2} \mid l$, making the right hand side of \eqref{ABtypeB} potentially nonzero. Suppose there exist $n_i$ such that $\frac{d}{2} \nmid n_i$. The first factor of the same expression can be decomposed as a product of $q$-binomial coefficients of the following form:
\[
\prod_{s = 1}^k \begin{bmatrix} \sum_{t = 1}^s n_t \\ n_s \end{bmatrix}_{q^2} .
\]
Up to permutation of the $n_p$'s, we are able to get the following $q$-binomial as a factor:
\begin{equation}
\label{rightbin}
\begin{bmatrix} l \\ n_i \end{bmatrix}_{q^2}.
\end{equation}
Using $q$-Lucas theorem to evaluate \eqref{rightbin} in $\zeta_d$ gives us a factor of the form
\[
\begin{bmatrix} 0 \\ s \end{bmatrix}_{q^2},
\]
where $s = n_i \mod \frac{d}{2}$, which is nonzero, then this factor is zero and so is the whole evaluation.
\item Suppose now that $\frac{d}{2} \mid l$ and $\frac{d}{2} \mid n_i$ for all $n_i$, let us show that these conditions imply $d \mid 2j+2$. Remember that $\rk(W_X) = n - l \equiv \frac{d}{2} - 1 \mod \frac{d}{2}$. We get that $\sum n_p i_p + j \equiv \frac{d}{2} - 1 \mod \frac{d}{2}$, which means that $j \equiv \frac{d}{2} - 1 \mod \frac{d}{2}$, so we get
$2j \equiv d - 2 \mod d$ and $d \mid 2j + 2$. In any other case $\mu_\lambda(\zeta_d) = 0$.
\end{enumerate}
\end{pr}

\begin{Lm}
\label{nonzeroB}
If $\frac{d}{2} \mid l$ and $\frac{d}{2} \mid n_i$ for every $n_i$, then: 
\[
\mu_\lambda(\zeta_d) = \binom{\frac{2l}{d}}{\frac{2n_1}{d}, \frac{2n_2}{d}, \cdots, \frac{2n_k}{d}} \binom{\frac{2n + 2}{d} + \frac{2l}{d}-1}{\frac{2l}{d}}.
\]
\end{Lm}
\begin{pr}
We will evaluate the first and second factor of \eqref{ABtypeB} separately:
\begin{itemize}
\item The first factor of \eqref{ABtypeB} can be written in the following form:
\[
\prod_{s = 1}^k \begin{bmatrix} \sum_{t = 1}^s n_t \\ n_s \end{bmatrix}_{q^2} .
\]
Using $q$-Lucas theorem on each of its subfactors, and knowing that $d$ divides each $n_i$ and $l$, we evaluated at $q = \zeta_d$, we get:
\[
\displaystyle \prod_{s = 1}^k \begin{bmatrix} 0 \\ 0 \end{bmatrix}_{\zeta_d^2} \binom{ \sum_{t = 1}^s \frac{2n_t}{d}}{\frac{2n_s}{d}},
\]
which simplifies as the following:
\[
\displaystyle \prod_{s = 1}^k \binom{ \sum_{t = 1}^s \frac{2n_t}{d}}{\frac{2n_s}{d}}.
\]
We can once again rewrite this expression as a multinomial coefficient:
\begin{equation}
\label{leftB}
\binom{ \frac{2l}{d}}{ \frac{2n_1}{d}, \frac{2n_2}{d}, \cdots, \frac{2n_k}{d} }.
\end{equation}
\item We will evaluate the second factor of \eqref{ABtypeB} at $q = \zeta_d$ using $q$-Lucas theorem, we get:
\[
\begin{bmatrix} \frac{d}{2} - 1 \\ \frac{d}{2} - 1 \end{bmatrix}_{\zeta_d^2} \binom{\frac{2n + 2}{d} + \frac{2l}{d}-1}{\frac{2n + 2}{d} - 1},
\]
which by the symmetry argument $\binom{n+k}{k} = \binom{n+k}{n}$ can be written as:
\[
\begin{bmatrix} \frac{d}{2} - 1 \\ \frac{d}{2} - 1 \end{bmatrix}_{\zeta_d^2} \binom{\frac{2n + 2}{d} + \frac{2l}{d}-1}{\frac{2l}{d}}.
\]
The remaining $q$-binomial is the unit polynomial, so we get
\begin{equation}
\label{rightB}
\binom{\frac{2n + 2}{d} + \frac{2l}{d}-1}{\frac{2l}{d}}.
\end{equation}
\end{itemize}
By combining \eqref{leftB} and \eqref{rightB}, we get the result.
\end{pr}

\subsection{Counting stable faces of $B_n$}

It now remains to count the faces which are stable under a $d$-fold rotation:
\begin{Lm}
\label{evenB2}
If $d$ is odd, then a type $B_n$ dissection is stable under a $d$-fold rotation if and only if it is stable under a $2d$-fold rotation.
\end{Lm}
\begin{pr}
The case $d = 1$ is trivial, so assume $d > 1$ is odd. If $f \in \Gamma(B_n)$ is stable under a $d$-fold rotation and under half-turn (which is automatic in type $B$), then it is stable under their composition $\rho$. It follows that $\rho$ is expressed as:
\[
e^{\frac{2i\pi}{d}} \cdot e^{i\pi} = e^{\frac{2(d+2)i\pi}{2d}}.
\]
As $d$ is odd, $d$ and $d+2$ are coprime, so $\rho$ has order $2d$.
\end{pr}
We can again restrict ourselves to the case where $d$ is even, and prove the "face counting version" of Lemmas \ref{zeroB} and \ref{nonzeroB}.
\begin{Lm}
\label{zeroB1}
If there exists $f \in \Gamma(B_n)_\lambda$ stable under $d$-fold rotation, then:
\begin{enumerate}
\item $\frac{d}{2} \mid l$,
\item $d \mid n_i$ for every $n_i$,
\item $d \mid 2j+2$.
\end{enumerate}
\end{Lm}
\begin{pr}
The proof is really straightforward. We only have to translate the propositions into the shape of the dissection:
\begin{enumerate}
\item $l$ is the number of subgons on each half of the dissection, except for the central polygon. If $\frac{d}{2} \nmid l$, a $d$-fold rotation cannot always send a polygon onto another one. As a consequence, the dissection cannot be stable under such a rotation.
\item Assume $\frac{d}{2} \mid l$, $n_i$ is the number of $(i + 2)$-gons on each half of the dissection, so if $\frac{d}{2} \nmid n_i$ for any $i$, a $d$-fold rotation cannot always send an $(i +2)$-gon onto another one. As a consequence, the dissection cannot be stable under such a rotation.
\item This is straightforward, the $B_j$ component of $W_X$ gives rise to a central $(2j+2)$-gon in the dissections, if $d \nmid 2j+2$, then a $d$-fold rotation cannot preserve the edges of this polygon.
\end{enumerate}
\end{pr}

In every other case, the number of faces stable under $d$-fold rotation, where $d$ is even, can be deduced from Adams and Banaian's formula in \cite{AB25}:
\begin{Lm}
\label{countB}
If $d$ is even, the number of faces of $\Gamma^{(m)}(B_n)_\lambda$ stable under $d$-fold rotation is given by:
\[
\binom{\frac{2l}{d}}{\frac{2n_1}{d}, \frac{2n_2}{d}, \cdots, \frac{2n_k}{d}} \binom{\frac{2n + 2}{d} + \frac{2l}{d}-1}{\frac{2l}{d}}.
\]
\end{Lm}

\begin{Thm}
The triple $(\Gamma^{(m)}(B_n)_\lambda, C_{2mn+2}, \mu_\lambda(q))$ exhibits the cyclic sieving phenomenon.
\end{Thm}
\begin{pr}
Lemmas \ref{evenB1} and \ref{evenB2} allow us restrict ourselves to the case where $d$ is even. By Lemmas \ref{nonzeroB}, \ref{countB}, \ref{zeroB} and \ref{zeroB1} the number of faces  in $\Gamma^{(m)}(B_n)_X$ stable under $d$-fold rotation is the same as $\mu_\lambda(\zeta_d)$ whenever $d$ is even. So we conclude the CSP.
\end{pr}

\section{Proof of the main theorem in dihedral types}
\label{PI}

\subsection{Combinatorics of the dihedral cluster complexes}

The complex $\Gamma^{(m)}(I_2(k))$ can be explicitly realized as the following complex of dissections of the $(mk+2)$-gon:
\begin{itemize}
\item a vertex of this complex is a family of parallel $m$-allowable diagonals such that the subgon defined by two such consecutive parallel diagonals is a $(2m+2)$-gon,
\item two vertices are \textit{compatible} if and only if all the diagonals of the dissection are noncrossing and form an $(m+2)$-angulation of the polygon.
\end{itemize}

\begin{figure}[h]
\begin{tikzpicture}[scale=2]
  \def\n{17}
  \def\r{1}
  \foreach \i in {1,...,\n} {
    \coordinate (P\i) at ({\r*cos(90 - 360/\n*(\i-1))}, {\r*sin(90 - 360/\n*(\i-1))});
  }
  \foreach \i in {1,...,\n} {
    \pgfmathtruncatemacro{\j}{mod(\i,\n)+1}
    \draw[thick] (P\i) -- (P\j);
  }
  \draw (P1) -- (P14);
  \draw (P4) -- (P11);
  \foreach \i in {1,...,\n} {
    \fill (P\i) circle (0.02);
  }
\end{tikzpicture}
    \begin{tikzpicture}[scale=2]
  \def\n{17}
  \def\r{1}
  \foreach \i in {1,...,\n} {
    \coordinate (P\i) at ({\r*cos(90 - 360/\n*(\i-1))}, {\r*sin(90 - 360/\n*(\i-1))});
  }
  \foreach \i in {1,...,\n} {
    \pgfmathtruncatemacro{\j}{mod(\i,\n)+1}
    \draw[thick] (P\i) -- (P\j);
  }
  \draw (P1) -- (P14);
  \draw (P4) -- (P11);
  \draw (P13) -- (P3);
  \draw (P6) -- (P10);
  \foreach \i in {1,...,\n} {
    \fill (P\i) circle (0.02);
  }
\end{tikzpicture}
\caption{A vertex and a facet of $\Gamma^{(3)}(I_2(5))$.}
\end{figure}
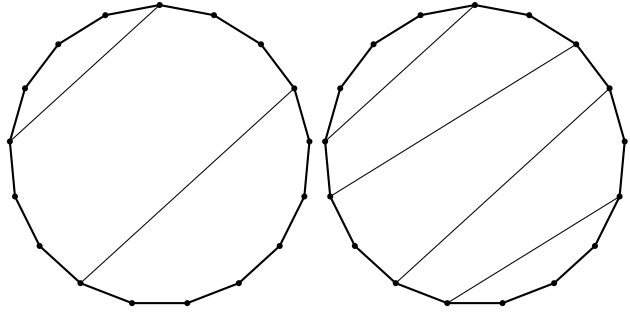

This combinatorial model relies on the inclusion $I_2(k) \hookrightarrow \mathfrak{S}_k$, see [\cite{DJV25}, Section 2.5] for more details.

\subsection{Sieving polynomials in dihedral types}

Let us compute our polynomial $\mu_\lambda$ in every case:
\[
\mu_\lambda(q) =
\begin{cases}
\frac{[1]_q}{[1]_q} &\text{ if $W_X = W$}, \\
\frac{[mk+2]_q}{[2]_q} &\text{ if $k$ is even and and $\rk(W_X) = 1$},\\
\frac{[mk+2]_q}{[1]_q} &\text{ if $k$ is odd and $\rk(W_X) = 1$},\\
\frac{[mk+2]_q [(m+1)k]_q}{[2]_q [k]_q} &\text{ if $W_X$ is trivial}.
\end{cases}
\]
\begin{rqe}
A cardinality argument shows that the two rank 1 standard parabolic subgroups of $W$ are parabolically conjugate if and only if $k$ is odd: there are exactly $k+2$ vertices in $\Gamma(I_2(k))$, which are split into two parabolic types if and only if $k$ is even.
\end{rqe}

\begin{Lm}
\label{I1}
If $d \mid mk+2$, $d > 2$ and $W_X \neq W$, then $\mu_\lambda(\zeta_d) = 0$.
\end{Lm}

\begin{pr} Let us split the cases by the rank of $W_X$:
\begin{itemize}
\item If $\rk(W_X) = 1$, as $d > 2$ the denominator is nonzero in any case. Then as $d \mid mk+2$, the numerator vanishes and so does $\mu_\lambda(\zeta_d)$.
\item If $\rk(W_X) = 0$, then the numerator vanishes as $d \mid mk+2$. Again as $d > 2$, we have $[2]_{q = \zeta_d} \neq 0$, and as $d \mid mk+2$, we have $mk \equiv d - 2 \mod d$, so that $d \nmid k$, and $[k]_{q = \zeta_d} \neq 0$. So $\mu_\lambda(\zeta_d) = 0$.
\end{itemize}
\end{pr}

\begin{Lm}
\label{I2}
If $k$ is odd and $mk+2$ is even, then:
\begin{itemize}
\item $\mu_\lambda(-1) = 0$ if $\rk(W_X) = 1$,
\item $\mu_\lambda(-1) = \frac{mk+2}{2}$ if $W_X$ is trivial,
\end{itemize}
and if $k$ is even:
\begin{itemize}
\item $\mu_\lambda(-1) = \frac{mk+2}{2}$ if $\rk(W_X) = 1$,
\item $\mu_\lambda(-1) = \frac{(mk+2)(m+1)}{2}$ if $W_X$ is trivial.
\end{itemize}
\end{Lm}

\begin{pr} Let us split the cases by rank again:
\begin{itemize}
\item Let $\rk(W_X) = 1$, as $mk+2$ is even we have $[mk+2]_{q=-1} = 0$, which concludes the case when $k$ is odd. When $k$ is even, as $mk+2$ and 2 are both even, then by \ref{q-prop} $\frac{[mk+2]_{q=-1}}{[2]_{q=-1}} = \frac{mk+2}{2}$.
\item Let $W_X$ be trivial, if we use the factorization $\mu_\lambda(q) = \frac{[mk+2]_q}{[2]_q} \frac{[(m+1)k]_q}{[k]_q}$, we get:
\begin{itemize}
\item By a preceeding argument, $\frac{[mk+2]_{q=-1}}{[2]_{q=-1}} = \frac{mk+2}{2}$.
\item If $k$ is even then so is $(m+1)k$, and then $\frac{[(m+1)k]_{q=-1}}{[k]_{q=-1}} = \frac{(m+1)k}{k}$.
\item If $k$ is odd and $mk+2$ is even, then $m$ is even too and therefore $m+1$ is odd. Then we get that $(m+1)k$ is odd and so $\frac{[(m+1)k]_{q=-1}}{[k]_{q=-1}} = 1$.
\end{itemize}
The result follows.
\end{itemize}
\end{pr}

\subsection{Counting stable faces in dihedral types}
In this section we consider $\Gamma^{(m)} = \Gamma^{(m)}(I_2(k))$ and $d \mid mk+2$. Let us count the faces of $\Gamma^{(m)}$ which are stable under $d$-fold rotation:

\begin{Lm}
\label{I3}
If $W_X \neq W$ and $d > 2$, then no $f \in \Gamma^{(m)}_X$ is stable under $d$-fold rotation.
\end{Lm}

\begin{pr}
This simply follows from proposition 1 of \cite{AB25}: such a dissection does not have a central polygon linked to more than 2 other subgons in its adjacency graph. Therefore there can be no $d$-fold rotation for $d > 2$.
\end{pr}

\begin{Lm} \label{I4} Let $\rk(W_X) = 1$:
\begin{itemize}
\item if $k$ is odd, then no $f \in \Gamma^{(m)}_X$ is stable under half-turn;
\item if $k$ is even, then every $f \in \Gamma^{(m)}_X$ is stable under half-turn.
\end{itemize}
\end{Lm}

\begin{pr}
We consider dissections made using only parallel diagonals. Let us consider the adjacency graph of the dissection: if $k$ is odd, its adjacency graph must have the following shape
\[
\begin{tikzpicture}[scale=1, every node/.style={draw=none}]
    \node (v1) {$m+2$};
    \node[right=of v1] (v2) {$2m+2$};
    \node[right=of v2, draw=none, shape=rectangle, minimum size=0mm] (v4) {$\dots$};
    \node[right=of v4] (v5) {$2m+2$};
    \node[right=of v5] (v7) {$2m+2$};
    \foreach \i/\j in {1/2,2/4,4/5,5/7}
        \draw (v\i) -- (v\j);
\end{tikzpicture}
\]
which cannot give rise to a half-turn as we would need the only $(m+2)$-gon to be fixed.\\
If $k$ is even, a vertex of $\Gamma^{(m)}$ can have either one of the following shapes:
\begin{center}
\begin{tikzpicture}[scale=1, every node/.style={draw=none}]
    \node (v1) {$m+2$};
    \node[right=of v1] (v2) {$2m+2$};
    \node[right=of v2, draw=none, shape=rectangle, minimum size=0mm] (v4) {$\dots$};
    \node[right=of v4] (v5) {$2m+2$};
    \node[right=of v5] (v7) {$m+2$};
    \foreach \i/\j in {1/2,2/4,4/5,5/7}
        \draw (v\i) -- (v\j);
\end{tikzpicture},
\end{center}
\begin{center}
\begin{tikzpicture}[scale=0.7, every node/.style={draw=none}]
    \node (v2) {$2m+2$};
    \node[right=of v2] (v3) {$2m+2$};
    \node[right=of v3, draw=none, shape=rectangle, minimum size=0mm] (v4) {$\dots$};
    \node[right=of v4] (v5) {$2m+2$};
    \node[right=of v5] (v6) {$2m+2$};
    \foreach \i/\j in {2/3,3/4,4/5,5/6}
        \draw (v\i) -- (v\j);
\end{tikzpicture}.
\end{center}

Such a set of diagonals is preserved by a half-turn, which acts nontrivially on the adjacency graph of the dissection. We get the result.
\end{pr}

\begin{Lm} \label{I5}
If $W_X$ is trivial, then:
\begin{itemize}
\item if $k$ is even, all facets of $\Gamma^{(m)}$ are stable under half-turn;
\item if $k$ is odd, exactly $\frac{mk+2}{2}$ facets of $\Gamma^{(m)}$ are stable under half-turn.
\end{itemize}
\end{Lm}

\begin{pr}
When $k$ is even, by the previous lemma all the vertices are preserved by the half-turn. Then any set of vertices must be preserved by half-turn too.\\
Now when $k$ is odd, let us notice that a $(m+2)$-angulation consisting only of parallel diagonals is preserved by such a rotation. Let us notice that as no vertex is stable under half-turn in that case, $(\alpha, \beta)$ is stable under half-turn if and only if $\alpha \mapsto \beta$ and $\beta \mapsto \alpha$. As the half-turn maps a diagonal into a parallel diagonal, a dissection stable by half-turn always has the described structure. The choice of such a dissection is equivalent to that of a class of parallel diagonals, which is itself equivalent to the choice of a diameter of the polygon. Therefore there are $\frac{mk+2}{2}$ such dissections.
\end{pr}

\begin{Thm}
The triple $(\Gamma^{(m)}(I_2(k))_X, C_{k+2}, \mu_X(q))$ exhibits the cyclic sieving phenomenon.
\end{Thm}

\begin{pr}
Combining Lemmas \ref{I1}, \ref{I2}, \ref{I3}, \ref{I4} and \ref{I5} gives us that the numbers we are comparing are the same. We conclude the CSP.
\end{pr}

\section{Proof of the main theorem in type $D_n$}
\label{PD}
In this section, let $W = D_n$ with Coxeter number $h = 2n-2$. We will again work directly on the generalized cluster complex. It is interesting to notice that for the first time we will not treat all parabolic conjugacy classes of parabolic subgroups. 
\subsection{Combinatorics of the type $D_n$ cluster complexes}
The combinatorics of the type $D_n$ objects is often a lot more irregular than the other classical types, and the cluster complexes reflect this. We follow the combinatorial interpretation of $\Gamma^{(m)}(D_n)$ given in \cite{FR05}. Let us consider the $2m(n-1)+2$-gon with counterclockwise labeled vertices $\{1, \cdots m(n-1)+1, \bar{1}, \cdots, \overline{m(n-1)+1} \}$. We call the edges $\{(mj+1, mj+2) \mid 0 \leq j \leq n-1  \}$ \textit{color-switching}.
\begin{itemize}
\item A \textit{type $D_n$ diagonal} is either a centrally symmetric pair of diagonals, or a diameter endowed with a color among two.
\item The Fomin-Reading rotation of the complex is the clockwise rotation of the polygon, sending the vertex $2$ on $1$, which changes the color of a diameter if it passes through a color-switching edge.
\item Two type $D_n$ diameters are compatible if and only if one of the following conditions hold:
\begin{enumerate}
\item One of them is a non-diameter diagonal and they are noncrossing.
\item They are confounded diameters of different color.
\item They are both non-confounded diameters and, up to rotation, when one of them is the diameter $\{1, -1 \}$ they have the same color.
\end{enumerate}
\end{itemize}
Classifying the faces $f \in \Gamma^{(m)}(D_n)$ by parabolic type is more complicated than in the previous types. We describe an algorithm.\\
Consider $\mu = [ \text{ } ]$ for the moment. For every subgon $P$ of $f$:
\begin{itemize}
\item If $P$ is an $(mk+2)$-gon which is not adjacent to exactly unicolored diameter, add a part $k$ to $\mu$.
\item If $P$ is an $(mk+2)$-gon, if there is a unique unicolored diameter in $f$ and if $P$ is adjacent to this diameter, add a part $(k+1)$ to $\mu$.
\item If $P$ is adjacent to two diameters it can be any $(mk+2+ i)$-gon, where $1 \leq i < m$. In that case, add a part $(k+1)$ to $\mu$.
\end{itemize}
At this stage, you get $\mu = [1^{n_1}, \cdots, k^{n_k}]$. The partition you get is $\lambda = [1^{\lfloor \frac{n_1}{2} \rfloor}, \cdots, k^{\lfloor \frac{n_k}{2} \rfloor}]$. If $\lambda \vdash n$ and $\lambda$ is all-even, add the sign $\pm$ corresponding to the $\mathcal{R}_m$-orbit of the diameter.

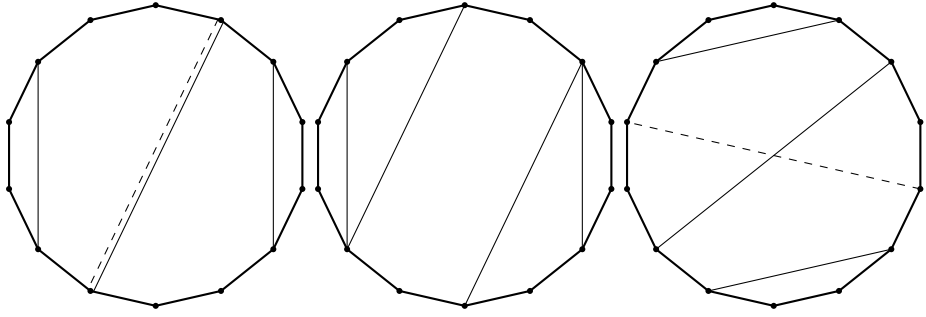
\begin{figure}[htbp]
\begin{tikzpicture}[scale=2]
  \def\n{14}
  \def\r{1}
  \foreach \i in {1,...,\n} {
    \coordinate (P\i) at ({\r*cos(90 - 360/\n*(\i-1))}, {\r*sin(90 - 360/\n*(\i-1))});
  }
  \foreach \i in {1,...,\n} {
    \pgfmathtruncatemacro{\j}{mod(\i,\n)+1}
    \draw[thick] (P\i) -- (P\j);
  }
\draw ($ (P2)!0.02cm!90:(P9) $) -- ($ (P9)!0.02cm!-90:(P2) $);
\draw[dashed] ($ (P2)!0.02cm!-90:(P9) $) -- ($ (P9)!0.02cm!90:(P2) $);
  \draw (P3) -- (P6);
  \draw (P10) -- (P13);
  \foreach \i in {1,...,\n} {
    \fill (P\i) circle (0.02);
  }
\end{tikzpicture}
\begin{tikzpicture}[scale=2]
  \def\n{14}
  \def\r{1}
  \foreach \i in {1,...,\n} {
    \coordinate (P\i) at ({\r*cos(90 - 360/\n*(\i-1))}, {\r*sin(90 - 360/\n*(\i-1))});
  }
  \foreach \i in {1,...,\n} {
    \pgfmathtruncatemacro{\j}{mod(\i,\n)+1}
    \draw[thick] (P\i) -- (P\j);
  }
  \draw (P1) -- (P10);
  \draw (P3) -- (P8);
  \draw (P3) -- (P6);
  \draw (P10) -- (P13);
  \foreach \i in {1,...,\n} {
    \fill (P\i) circle (0.02);
  }
\end{tikzpicture}
\begin{tikzpicture}[scale=2]
  \def\n{14}
  \def\r{1}
  \foreach \i in {1,...,\n} {
    \coordinate (P\i) at ({\r*cos(90 - 360/\n*(\i-1))}, {\r*sin(90 - 360/\n*(\i-1))});
  }
  \foreach \i in {1,...,\n} {
    \pgfmathtruncatemacro{\j}{mod(\i,\n)+1}
    \draw[thick] (P\i) -- (P\j);
  }
\draw[dashed] (P5) -- (P12);
\draw (P10) -- (P3);
\draw (P13) -- (P2);
\draw (P6) -- (P9);
  \foreach \i in {1,...,\n} {
    \fill (P\i) circle (0.02);
  }
\end{tikzpicture}
\caption{Faces of $\Gamma^{(2)}(D_3)$.}
\end{figure}

\subsection{Sieving polynomials in type $D_n$}
Unlike in the previous types, the sieving polynomials in type $D_n$ are not so uniform. Due to $N(W_X)/W_X$ not always being a reflection group (see Subsection \ref{Ddenom}), we lack a $q$-analogue of our formula in some families. We treat the following cases seperately:
\begin{enumerate}
\item \label{1} $\lambda \vdash k < n$;
\item \label{2} $(\lambda, \pm)$, $\lambda \vdash n$ all even;
\item \label{3} $\lambda = [2^{n_2}, 3^{n_3}, \cdots, k^{n_k}] \vdash n$, with unique odd $k \ne 1$ such that $n_m > 0$ ;
\item \label{4} $\lambda = [2^{n_2}, 3^{n_3}, \cdots, k^{n_k}] \vdash n$, with $n_m = 0$ for every odd $m \ne 1$.
\end{enumerate}

Let us first treat case \ref{1}, namely $\lambda \vdash k < n$. Let us recall that we consider $\rk(W_X) = n - l$. Section \ref{q-counting} and Subsection \ref{Ddenom} show us that we have $\Pi(D_n)_X \cong \Pi(B_{\dim X})$ and $N(W_X)/W_X \cong \prod_{i \in \lambda} B_{n_i}$. The polynomial $\mu_\lambda$ in case \ref{1} can be computed as follows:
\begin{align*}
\displaystyle \mu_\lambda(q) &= \frac{\prod_{i=1}^{l}[2m(n-1)+2i]_q}{[2n_1]_q!! \cdots [2n_k]_q!!},\\
&= \frac{\prod_{i=1}^{l}[m(n-1)+i]_{q^2}}{[n_1]_{q^2}! \cdots [n_k]_{q^2}!},\\
&= \frac{[m(n-1)+l]_{q^2}!}{[n_1]_{q^2}! \cdots [n_k]_{q^2}! [m(n-1)]_{q^2}!},\\
&= \begin{bmatrix} m(n-1) +l \\ n_1, \cdots, n_k, m(n-1) \end{bmatrix}_{q^2}.
\end{align*}
This polynomial is exactly the sieving polynomial $\mu_\lambda$ for faces of $\Gamma^{(m)}(B_{n-1})$ of type $\lambda$, thus its evaluation on the roots of unity of interest for us is detailed in Subsection \ref{polyB}.\\
Case \ref{2} is also quite easy. Let $\lambda \vdash n$ be an all-even partition. Section \ref{q-counting} and Subsection \ref{Ddenom} show us again that $\Pi(D_n)_X \cong \Pi(B_{\dim X})$ and $N(W_X)/W_X \cong \prod_{i \in \lambda} B_{n_i}$. Then again we get
\[
\mu_\lambda(q) = \begin{bmatrix} m(n-1) +l \\ n_2, n_4, \cdots, n_k, m(n-1) \end{bmatrix}_{q^2}
\]
which evaluates in the same way as the preceding case by Subsection \ref{polyB}.\\
Case \ref{3} is given by $\lambda \vdash n$ with only one nonzero $n_i$ for $i$ odd $\neq 1$. Let $\lambda$ be such a partition, and fix $k$ the odd integer such that $n_i > 0$. Section \ref{q-counting} and Subsection \ref{Ddenom} show us that $\Pi(D_n)_X \cong \Pi(B_{\dim X})$ and $N(W_X)/W_X \cong D_{n_k} \times \prod_{i \neq k} B_{n_i}$. In that case $\mu_X$ is given by
\[
\displaystyle \mu_\lambda(q) = \frac{\prod_{i=1}^{l}[2m(n-1)+2i]_q}{[2n_2]_q!! [2n_4]_q!! \cdots [2n_k-1]_q!! [n_k]_q}.
\]
Multiply by $\frac{[2n_k]_q}{[2n_k]_q}$ and rearrange to get 
\begin{equation}
\label{polyD3}
\displaystyle \mu_\lambda(q) = \begin{bmatrix} m(n-1) + l \\ n_2, n_4, \cdots, n_k, m(n-1) \end{bmatrix}_{q^2} \frac{[2n_k]_q}{[n_k]_q}.
\end{equation}
\begin{Lm}
\label{oddpartneq1}
If $\lambda \vdash n$ belongs in case \ref{3} and $d \mid 2m(n-1)+2$, then:
\[
\mu_\lambda(\zeta_d) = \begin{cases} 2 \binom{\frac{m(n-1)+1}{d}+\frac{l}{d} - 1}{\frac{n_2}{d}, \frac{n_4}{d}, \cdots, \frac{n_k}{d}, \frac{m(n-1)+1}{d}-1} &\text{ if } d \mid n_i \text{ for every } n_i \text{ and } d \text{ is odd},\\
2 \binom{\frac{2m(n-1)+2}{d}+\frac{2l}{d} - 1}{\frac{2n_2}{d}, \frac{2n_4}{d}, \cdots, \frac{2n_k}{d}, \frac{2m(n-1)+2}{d}-1} &\text{ if } d \text{ is even, } \frac{d}{2} \mid n_i \text{ for every } n_i \text{ and } d \mid n_k,\\
0 &\text{ otherwise.}\end{cases}
\]
\end{Lm}
\begin{pr}
The left hand side of expression \eqref{polyD3} can be computed using Lemmas 6.1, 6.2 and by replacing $n$ by $m(n-1)$. By splitting cases according to the parity of $d$, we either find that if the left hand side of \eqref{polyD3} is nonzero when evaluated in $\zeta_d$, then $d \mid m(n-1) + 1$, $ d \mid l$ and $d \mid n_i$ for every $n_i$; or we find that $\frac{d}{2} \mid m(n-1) + 1$, $\frac{d}{2} \mid l$ and $\frac{d}{2} \mid n_i$ for every $n_i$. When $d$ is odd, then the right hand side of \eqref{polyD3} automatically evaluates to 2 whenever the left hand side is nonzero, thanks to Proposition \ref{q-prop} and the fact that $d \mid n_k$. If $d$ is even, the same argument applies if $d \mid n_k$. However if $d \nmid n_k$, then $\frac{[2n_k]_q}{[n_k]_q}$ evaluates to 0 in $\zeta_d$.
\end{pr}
Let us finally treat case \ref{4}, given by $\lambda \vdash n$ where $n_k = 0$ for odd $k \neq 1$. Section \ref{q-counting} and Subsection \ref{Ddenom} show us that we have $\Pi(D_n)_X \cong \Pi(D_{\dim X}^{l(\lambda) - n_1})$ and $N(W_X)/W_X \cong D_{n_1} \times \prod_{i \geq 2} B_{n_i}$. The polynomial $\mu_\lambda$ in that case can be computed as follows:
\[
\displaystyle \mu_\lambda(q) = \frac{\prod_{i=1}^{l-1}[2m(n-1)+2i]_q}{[2n_1-2]_q!! \cdots [2n_k]_q!!} \frac{[2m(n-1)+ 2l - n_1]_q}{[n_1]_q}.
\]
Which can be seen as:
\begin{equation}
\label{D11}
\mu_\lambda(q) = \begin{bmatrix} m(n-1) +l -1 \\ n_1-1, \cdots, n_k, m(n-1) \end{bmatrix}_{q^2} \frac{[2m(n-1)+ 2l - n_1]_q}{[n_1]_q}.
\end{equation}
By adding factors $\frac{[2m(n-1) + 2l]_q}{[2m(n-1) + 2l]_q}$ and $\frac{[2n_1]_q}{[2n_1]_q}$, we get
\begin{equation}
\label{D12}
\mu_\lambda(q) = \begin{bmatrix} m(n-1) +l \\ n_1, \cdots, n_k, m(n-1) \end{bmatrix}_{q^2} \frac{[2m(n-1)+ 2l - n_1]_q}{[2m(n-1) + 2l]_q} \frac{[2n_1]_q}{[n_1]_q}.
\end{equation}

\begin{Lm}
\label{oddpart1}
If $\lambda$ belongs in case \ref{4}, then:
\begin{itemize}
\item if $n$ is odd, then $\mu_\lambda(-1) = \binom{ m(n-1) +l -1 }{ n_1-1, n_2 , n_4, \cdots, m(n-1) } $;
\item if $n$ is even, then $\mu_\lambda(-1) = \mu_\lambda(1)$;
\item if $d > 2$, $d \mid n$ and $d \mid n_i$ for every $i \in \lambda$, then
\[
\mu_\lambda(\zeta_d) = \begin{cases} 2 \binom{\frac{m(n-1)+1}{d}+\frac{l}{d} - 1}{\frac{n_1}{d}, \frac{n_2}{d}, \frac{n_4}{d}, \cdots, \frac{m(n-1)+1}{d}-1} &\text{ if } d \mid n_i \text{ for every } n_i \text{ and } d \text{ is odd},\\
2 \binom{\frac{2m(n-1)+2}{d}+\frac{2l}{d} - 1}{ \frac{2n_1}{d}, \frac{2n_2}{d}, \frac{2n_4}{d}, \cdots, \frac{2m(n-1)+2}{d}-1} &\text{ if } d \text{ is even, } \frac{d}{2} \mid n_i \text{ for every } n_i \text{ and } d \mid n_1,\\
0 &\text{ otherwise.}\end{cases}
\]
\end{itemize}
\end{Lm}

\begin{pr}Let us first treat the case $q = -1$, and then $q = \zeta_d$ with $d > 2$.
\begin{itemize}
\item Let us notice that $n$ and $n_1$ have the same parity. If $n$ is even, it is therefore clear that all the $q$-integers involved in $\mu_\lambda$ are even and then it follows easily that $\mu_\lambda(1) = \mu_\lambda(-1)$. If $n$ is odd, then so is $n_1$ and we get that $2m(n-1) + 2l - n_1$ and $n_1$ are both odd, meaning that $\frac{[2m(n-1)+ 2l - n_1]_{q=-1}}{[n_1]_{q=-1}} = 1$. It follows in that case using expression \eqref{D11} that $\mu_\lambda(-1) = \binom{ m(n-1) +l -1 }{ n_1-1, n_2, n_4 \cdots, m(n-1) }$.
\item Suppose now $d > 2$. Using expression \eqref{D12} we already know by the previous cases how to evaluate
\[
\begin{bmatrix} m(n-1) +l \\ n_1, \cdots, n_k, m(n-1) \end{bmatrix}_{q^2} \frac{[2n_1]_q}{[n_1]_q}
\]
in $q = \zeta_d$. It remains to evaluate
\[
\frac{[2m(n-1)+ 2l - n_1]_q}{[2m(n-1) + 2l]_q}
\]
in $q = \zeta_d$. This is done quickly by noticing that as $d \mid n_1$ if $\mu_X(\zeta_d) \neq 0$, $2m(n-1)+ 2l - n_1 \equiv 2m(n-1) + 2l \equiv d - 2 \mod d$. It follows that 
\[
\frac{[2m(n-1)+ 2l - n_1]_{q= \zeta_d}}{[2m(n-1) + 2l]_{q = \zeta_d}} = 1.
\]
We get the expected result.
\end{itemize}
\end{pr}

\subsection{Counting stable faces of $D_n$}
Counting the faces of $\Gamma^{(m)}(D_n)$ stable by a $d$-fold rotation is again much more complicated than the previous cases. The classification of the faces by parabolic type strictly separates the faces with at least one diameter (given by cases \ref{2}, \ref{3} and \ref{4}) and the ones without diameter (given by case \ref{1}). Case \ref{1} is treated by the following Lemma.

\begin{Lm}
\label{C1}
Let $\lambda \vdash m < n - 1$.
If $d$ is even, the faces of type $\lambda$ stable under $d$-fold rotation are counted by
\[
\binom{\frac{2l}{d}}{\frac{2n_1}{d}, \frac{2n_2}{d}, \cdots, \frac{2n_k}{d}} \binom{\frac{2m(n-1) + 2}{d} + \frac{2l}{d}-1}{\frac{2l}{d}}
\]
if $d \mid 2n_i$ for every $i$, and $0$ otherwise.\\
If $d$ is odd, then a face of type $\lambda$ is stable by a $d$-fold rotation if and only if it is stable by a $2d$-fold rotation.
\end{Lm}
\begin{pr}
Let us notice that the faces of type $\lambda$ can be described as the faces of the same type in $\Gamma^{(m)}(B_{n-1})$. The result follows from Lemmas \ref{evenB2} and \ref{countB}.
\end{pr}
We now have to deal with the existence of a diameter in the dissection, the two following Lemmas will be useful.
\begin{Lm}[\cite{EF08}, Lemma 5.6]
If $d \geq 3$, $d \mid 2m(n-1)+2$, $d \mid 2l$ and there exists $f \in \Gamma^{(m)}(D_n)$ of dimension $l-1$ stable under $d$-fold rotation with at least one diameter, then $d \mid n$.
\end{Lm}
\begin{rqe} Notice the following:
\begin{itemize}
\item If $f \in \Gamma^{(m)}(D_n)$ is stable under a $d$-fold rotation, then we must have $d \mid 2m(n-1)+2$ and $d \mid 2l$.
\item Given $d \mid 2m(n-1)+2$, the condition $d \mid n$ is equivalent to $\frac{2n}{d}$ being even and to $\lfloor \frac{2m(n-1)+2}{dm} \rfloor$ being even as well.
\end{itemize}
\end{rqe}

\begin{Lm}
\label{reduction}
Let $d \mid 2m(n-1)+2$ such that $d \geq 3$ and let $\lambda = \left[1^{n_1}, \cdots, k^{n_k}\right]$. Fix a color for the diameter $P_1 = \{-1,1\}$ and fix the color of the other diameters by taking the color that makes it compatible with $P_1$. The subset $G_{d, \lambda, \pm} \subset \Gamma^{(m)}(D_n)_\lambda$ formed by dissections with diameters of the right color that are stable under (purely geometrical) $d$-fold rotation are in bijection with:
\begin{enumerate}
\item the faces of $\Gamma^{(m)}\left(D_{\frac{n}{d}}\right)_{\lambda'}$ with at least a diameter of the right color and $\lambda'  = \left[1^{\frac{n_1}{d}}, \cdots, k^{\frac{n_k}{d}}\right]$ if $n$ is odd and $d \mid n_i$ for every $i$;
\item the faces of $ \displaystyle \Gamma^{(m)}\left(D_{\frac{2n}{d}}\right)_{\lambda'}$ with at least a diameter of the right color and $\displaystyle \lambda' = \left[1^{\frac{2n_1}{d}}, \cdots, k^{\frac{2n_k}{d}}\right]$ if $n$ is even and $d \mid 2n_i$ for every $i$;
\item $\emptyset$ otherwise.
\end{enumerate}
\end{Lm}

\begin{rqe}
We authorize a diameter to be bicolored in the image dissection, that may happen in the case where two diameters end up being confounded.
\end{rqe}

\begin{pr}
Let us look at a necessary condition for such a dissection to be stable under $d$-fold rotation:
\begin{itemize}
\item When $d$ is odd, because of the half-turn stability condition we need the diameters to divide the $(2m(n-1)+2)$-gon into $2d$ identical regions. This can be formed using $d$ diameters such that two consecutive of these diameters form a $\frac{\pi}{d}$ angle.
\item When $d$ is even, we need the diameters to divide the $(2m(n-1)+2)$-gon into $d$ identical regions. This can be formed using $\frac{d}{2}$ diameters such that two consecutive of these diameters form a $\frac{2\pi}{d}$ angle.
\end{itemize}
In every case we need the divisions to make sense because all the regions must be equivalent, so the classification of faces by their parabolic type forces these conditions to hold.
\begin{figure}[h!]
\begin{tikzpicture}[scale=2]
  \def\n{66}
  \def\r{1}
  \foreach \i in {1,...,\n} {
    \coordinate (P\i) at ({\r*cos(90 - 360/\n*(\i-1))}, {\r*sin(90 - 360/\n*(\i-1))});
  }
  \foreach \i in {1,...,\n} {
    \pgfmathtruncatemacro{\j}{mod(\i,\n)+1}
    \draw[thick] (P\i) -- (P\j);
  }
  \draw[dashed] (P2) -- (P35);
  \draw (P24) -- (P57);
  \draw[dashed] (P13) -- (P46);
  \draw (P4) -- (P9);
  \draw (P26) -- (P31);
  \draw (P37) -- (P42);
  \draw (P59) -- (P64);
  \draw (P15) -- (P20);
  \draw (P48) -- (P53);
  \foreach \i in {1,...,\n} {
    \fill (P\i) circle (0.02);
  }
\end{tikzpicture}
\caption{A face of $\Gamma^{(4)}(D_9)_{[1^3,2^3]}$ stable under $3$-fold rotation (not 6-fold because of the color-switches).}
\end{figure}
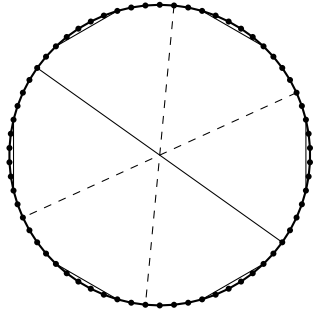

Consider the vertices of the $(2m(n-1)+2)$-gon to be roots of unity, the data of such a dissection is then considered to be the data of a small portion of the unit circle determined by a region. By applying the homeomorphism $z \mapsto z^d$ or $z \mapsto z^{\frac{d}{2}}$ depending on the parity of $d$, the families of dissections that we get are exactly the ones we wanted.\\
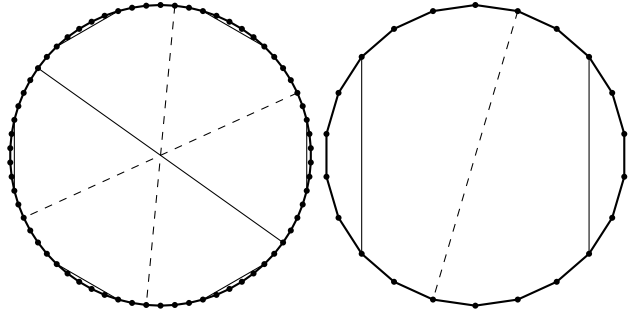
\begin{figure}[h!]
\begin{tikzpicture}[scale=2]
  \def\n{66}
  \def\r{1}
  \foreach \i in {1,...,\n} {
    \coordinate (P\i) at ({\r*cos(90 - 360/\n*(\i-1))}, {\r*sin(90 - 360/\n*(\i-1))});
  }
  \foreach \i in {1,...,\n} {
    \pgfmathtruncatemacro{\j}{mod(\i,\n)+1}
    \draw[thick] (P\i) -- (P\j);
  }
  \draw[dashed] (P2) -- (P35);
  \draw (P24) -- (P57);
  \draw[dashed] (P13) -- (P46);
  \draw (P4) -- (P9);
  \draw (P26) -- (P31);
  \draw (P37) -- (P42);
  \draw (P59) -- (P64);
  \draw (P15) -- (P20);
  \draw (P48) -- (P53);
  \foreach \i in {1,...,\n} {
    \fill (P\i) circle (0.02);
  }
\end{tikzpicture}
\begin{tikzpicture}[scale=2]
  \def\n{22}
  \def\r{1}
  \foreach \i in {1,...,\n} {
    \coordinate (P\i) at ({\r*cos(90 - 360/\n*(\i-1))}, {\r*sin(90 - 360/\n*(\i-1))});
  }
  \foreach \i in {1,...,\n} {
    \pgfmathtruncatemacro{\j}{mod(\i,\n)+1}
    \draw[thick] (P\i) -- (P\j);
  }
  \draw (P4) -- (P9);
  \draw (P15) -- (P20);
  \draw[dashed] (P2) -- (P13);
  \foreach \i in {1,...,\n} {
    \fill (P\i) circle (0.02);
  }
\end{tikzpicture}
\caption{A face of $\Gamma^{(4)}(D_9)_{[1^3,2^3]}$ stable under $3$-fold rotation and its image by the bijection.}
\end{figure}
\end{pr}

Using Lemma \ref{reduction}, we now treat case \ref{2}:

\begin{Cor}
\label{C2}
Let $d \mid 2m(n-1)+2$.
If $d$ is even, the set of faces of type $(\lambda, \pm)$, where $\lambda \vdash n$ is all-even, stable by a $d$-fold rotation is counted by:
\[
\binom{\frac{2l}{d}}{\frac{2n_2}{d}, \cdots, \frac{2n_k}{d}} \binom{ \frac{2m(n-1)+2}{d} + \frac{2l}{d} - 1}{\frac{2l}{d}}
\]
if $d \mid 2n_i$ for every $i$, and $0$ otherwise.\\
If $d$ is odd, then it is counted by:
\[
\binom{\frac{l}{d}}{\frac{n_2}{d}, \cdots, \frac{n_k}{d}} \binom{ \frac{m(n-1)+1}{d} + \frac{l}{d} - 1}{\frac{l}{d}}
\]
if $d \mid n_i$ for every $i$, and $0$ otherwise.
\end{Cor}

\begin{pr}
This follows straight from Lemma \ref{reduction}:\begin{itemize}
\item The set we are counting is exactly $G_{d, \lambda, \pm}$.
\item The constraints on $d$ and $\lambda$ force the condition $d \mid n$. Thus the $d$-fold Fomin-Reading rotation preserves the color of each diameter in every $f \in G_{d, \lambda, \pm}$.
\item The bijection from Lemma \ref{reduction} is in that case a bijection between $G_{d, \lambda, \pm}$ and $\Gamma^{(m)}(D_{\frac{2n}{d}})_{\lambda'}$ or $\Gamma^{(m)}(D_{\frac{2n}{d}})_{\lambda'}$ depending on the parity of $d$.
\end{itemize}
\end{pr}

We will also treat case \ref{3} using Lemma \autoref{reduction}:

\begin{Cor}
\label{C3}
Let $d \mid 2m(n-1)+2$.
If $d$ is even, the set of faces of type $\lambda$, where $\lambda$ belongs in case \ref{3}, stable by a $d$-fold rotation is counted by:
\[
2 \binom{\frac{2l}{d}}{\frac{2n_2}{d}, \cdots, \frac{2n_k}{d}} \binom{ \frac{2m(n-1)+2}{d} + \frac{2l}{d} - 1}{\frac{2l}{d}}
\]
if $d \mid 2n_i$ for every $i$, and $0$ otherwise.\\
If $d$ is odd, then it is counted by:
\[
2 \binom{\frac{l}{d}}{\frac{n_2}{d}, \cdots, \frac{n_k}{d}} \binom{ \frac{m(n-1)+1}{d} + \frac{l}{d} - 1}{\frac{l}{d}}
\]
if $d \mid n_i$ for every $i$, and $0$ otherwise.
\end{Cor}

\begin{pr}
This follows straight from Lemma \ref{reduction}:\begin{itemize}
\item The set we are counting is exactly $G_{d, \lambda, +} \uplus G_{d, \lambda, -}$.
\item The constraints on $d$ and $\lambda$ force the condition $d \mid n$. Thus the $d$-fold Fomin-Reading rotation preserves the color of each diameter in every $f \in G_{d, \lambda, \pm}$.
\item The bijection from Lemma \ref{reduction} is in that case a bijection between $G_{d, \lambda, +} \uplus G_{d, \lambda, -}$ and $\Gamma^{(m)}(D_{\frac{2n}{d}})_{\lambda'}$ or $\Gamma^{(m)}(D_{\frac{2n}{d}})_{\lambda'}$ depending on the parity of $d$.
\end{itemize}
\end{pr}

Case \ref{4} is again treated using Lemma \ref{reduction}:

\begin{Cor}
\label{C4}
Let $\lambda \vdash n$ such that $\lambda$ belongs in case \ref{4}, and let $d \mid 2m(n-1)+2$. The faces of $\Gamma^{(m)}(D_n)_\lambda$ stable by $d$-fold rotation are counted by:
\begin{itemize}
\item $\binom{ m(n-1) +l -1 }{ n_1-1, n_2 , n_4, \cdots, m(n-1) } $ if $n$ is odd and $d=2$;
\item $\mu_\lambda$ if $n$ is even and $d=2$;
\item $2 \binom{\frac{m(n-1)+1}{d}+\frac{l}{d} - 1}{\frac{n_1}{d}, \frac{n_2}{d}, \frac{n_4}{d}, \cdots, \frac{m(n-1)+1}{d}-1}$ if $d > 2$ is odd and $d \mid n_i$ for every $i \in \lambda$;
\item $2 \binom{\frac{2m(n-1)+2}{d}+\frac{2l}{d} - 1}{ \frac{2n_1}{d}, \frac{2n_2}{d}, \frac{2n_4}{d}, \cdots, \frac{2m(n-1)+2}{d}-1}$ if $d > 2$ is even and $\frac{d}{2} \mid n_i$ for every $i \in \lambda$ and $d \mid n_1$;
\item $0$ otherwise.
\end{itemize}
\end{Cor}

\begin{pr}
If $d = 2$, then this is quick: if $n$ is even then every diameter is stable under half-turn, and if $n$ is odd then no diameter is stable under half-turn. Thus for $n$ odd, the only possibility is to have a bicolored diameter. This is equivalent to counting classical dissections with a unique diameter, so by \cite{AB25} this set is counted by $\binom{ m(n-1) +l -1 }{ n_1-1, n_2 , n_4, \cdots, m(n-1) } $.
This follows straight from Lemma \ref{reduction}:\begin{itemize}
\item The set we are counting is exactly $G_{d, \lambda, +} \uplus G_{d, \lambda, -}$.
\item The constraints on $d$ and $\lambda$ force the condition $d \mid n$. Thus the $d$-fold Fomin-Reading rotation preserves the color of each diameter in every $f \in G_{d, \lambda, \pm}$.
\item The bijections from Lemma \ref{reduction} induce a surjection $G_{d, \lambda, +} \uplus G_{d, \lambda, -} \twoheadrightarrow \Gamma^{(m)}(D_{\frac{2n}{d}})_{\lambda'}$ or $G_{d, \lambda, +} \uplus G_{d, \lambda, -} \twoheadrightarrow \Gamma^{(m)}(D_{\frac{2n}{d}})_{\lambda'}$ depending on the parity of $d$. This function fails to be bijective because the image of a dissection with inner triangles formed with two diameters of the same $\mathcal{R}_m^{\frac{2m(n-1)+2}{d}}$-orbit or $\mathcal{R}_m^{\frac{m(n-1)+1}{d}}$-orbit (depending on the parity of $d$) does not depend on the choice of $\mathcal{R}_m$-orbit of diameter we choose for the colors. The images of these dissections happen to be the dissections of $\Gamma^{(m)}(D_{\frac{2n}{d}})_{\lambda'}$ or $\Gamma^{(m)}(D_{\frac{n}{d}})_{\lambda'}$ with a bicolored diameter, thus we have to count these specific dissections with multiplicity $2$. We know the cardinality $k$ of this set depending on the parity of $d$, and a quick calculation of the sum $| \Gamma^{(m)}(D_{\frac{2n}{d}})_{\lambda'}| + k$ gives us the result.
\end{itemize}
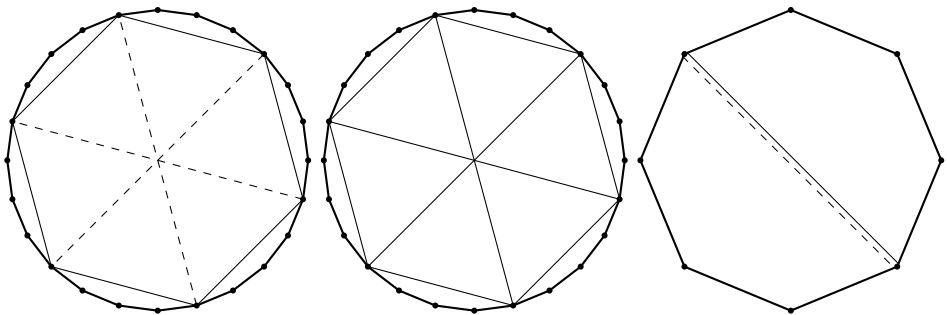
\begin{figure}[htbp]
\begin{tikzpicture}[scale=2]
  \def\n{24}
  \def\r{1}
  \foreach \i in {1,...,\n} {
    \coordinate (P\i) at ({\r*cos(90 - 360/\n*(\i-1))}, {\r*sin(90 - 360/\n*(\i-1))});
  }
  \foreach \i in {1,...,\n} {
    \pgfmathtruncatemacro{\j}{mod(\i,\n)+1}
    \draw[thick] (P\i) -- (P\j);
  }
  \draw[dashed] (P12) -- (P24);
  \draw[dashed] (P4) -- (P16);
  \draw[dashed] (P8) -- (P20);
  \draw (P4) -- (P8);
  \draw (P8) -- (P12);
  \draw (P12) -- (P16);
  \draw (P16) -- (P20);
  \draw (P20) -- (P24);
  \draw (P4) -- (P24);
  \foreach \i in {1,...,\n} {
    \fill (P\i) circle (0.02);
  }
\end{tikzpicture}
\begin{tikzpicture}[scale=2]
  \def\n{24}
  \def\r{1}
  \foreach \i in {1,...,\n} {
    \coordinate (P\i) at ({\r*cos(90 - 360/\n*(\i-1))}, {\r*sin(90 - 360/\n*(\i-1))});
  }
  \foreach \i in {1,...,\n} {
    \pgfmathtruncatemacro{\j}{mod(\i,\n)+1}
    \draw[thick] (P\i) -- (P\j);
  }
  \draw (P12) -- (P24);
  \draw (P4) -- (P16);
  \draw (P8) -- (P20);
  \draw (P4) -- (P8);
  \draw (P8) -- (P12);
  \draw (P12) -- (P16);
  \draw (P16) -- (P20);
  \draw (P20) -- (P24);
  \draw (P4) -- (P24);
  \foreach \i in {1,...,\n} {
    \fill (P\i) circle (0.02);
  }
\end{tikzpicture}
\begin{tikzpicture}[scale=2]
  \def\n{8}
  \def\r{1}
  \foreach \i in {1,...,\n} {
    \coordinate (P\i) at ({\r*cos(90 - 360/\n*(\i-1))}, {\r*sin(90 - 360/\n*(\i-1))});
  }
  \foreach \i in {1,...,\n} {
    \pgfmathtruncatemacro{\j}{mod(\i,\n)+1}
    \draw[thick] (P\i) -- (P\j);
  }
\draw ($ (P8)!0.02cm!90:(P4) $) -- ($ (P4)!0.02cm!-90:(P8) $);
\draw[dashed] ($ (P8)!0.02cm!-90:(P4) $) -- ($ (P4)!0.02cm!90:(P8) $);
  \foreach \i in {1,...,\n} {
    \fill (P\i) circle (0.02);
  }
\end{tikzpicture}
\caption{Two faces of $\Gamma^{(1)}(D_{12})$ on the left and the middle which both have the right dissection as an image by the bijections of Lemma \ref{reduction}.}
\end{figure}
\end{pr}

\begin{Thm}
Let $ \lambda \vdash m < n$ or $\lambda \vdash n$ with unique odd $i$ such that $n_i > 0$. The triple $(\Gamma^{(m)}(D_n)_\lambda, C_{2m(n-1)+2}, \mu_X(q))$ exhibits the cyclic sieving phenomenon.
\end{Thm}

\begin{pr}
Let us consider the cases separately:
\begin{enumerate}
\item For case \ref{1}, combine Lemmas \ref{evenB1}, \ref{zeroB}, \ref{nonzeroB} and \ref{C1}.
\item For case \ref{2}, combine Lemmas \ref{evenB1}, \ref{zeroB}, \ref{nonzeroB} and Corollary \ref{C2}.
\item For case \ref{3}, combine Lemma \ref{oddpartneq1} and Corollary \ref{C3}.
\item For case \ref{4}, combine Lemma \ref{oddpart1} and Corollary \ref{C4}.
\end{enumerate}
\end{pr}

\section{Proof of the main theorem in the exceptional types}
\label{exceptional}
At this point there only remains a finite number of classical irreducible cluster complexes, so we can easily computer-check our result in the case $m=1$. We use the tables in \cite{H80} and \cite{DPR25} as a reference for the structure of the groups $N(W_X)/W_X$ and \cite{OS83} for the exponents of $\Pi(W)_X$. In every table, we list the parabolic classes which satisfy the reflection property, except for the Catalan case which was proven in \cite{EF08} and the group $W$ itself which always yields $1$ (the empty face), and give the corresponding sieving polynomial. The reader who wishes to know the number of faces stable under $d$-fold rotation should use Proposition \ref{q-prop} to evaluate the polynomial in $\zeta_d$ in all the divisors of $h+2$.\\
Some groups in these tables do not satisfy the reflection property, but using the \cite{DPR25} decomposition, we have been able to find a suitable $q$-analogue by taking the reflection action on $X \cap Y$ and $Y^\perp$. These groups are preceeded in the tables by a *.

\begin{table}[!h]
\centering
\renewcommand{\arraystretch}{1.3}
\begin{tabular}{@{} l l c @{}}
\toprule
\textbf{$\lambda$} & \textbf{$\mu_\lambda(q)$} & \textbf{Divisors of $h+2$} \\
\midrule

$A_1$ & $\frac{[16]_q [12]_q}{[2]_q  [2]_q}$ & \multirow{2}{*}{\parbox[c]{3cm}{\centering
1,2,3,4,6,12}} \rule{0pt}{4ex}\\

$A_1^2, A_2, I_2(5)$ & $\frac{[12]_q}{ [2]_q}$ & \rule{0pt}{4ex}\\
\bottomrule
\end{tabular}
\caption{Sieving polynomials and divisors in type $H_3$}
\end{table}

\begin{table}[!h]
\centering
\renewcommand{\arraystretch}{1.3}
\begin{tabular}{@{} l l c @{}}
\toprule
\textbf{$\lambda$} & \textbf{$\mu_\lambda(q)$} & \textbf{Divisors of $h+2$} \\
\midrule

$A_1$ & $\frac{[50]_q [42]_q [32]_q}{[2]_q  [6]_q [10]_q}$ & \multirow{6}{*}{\parbox[c]{3cm}{\centering
1,2,4,8,16,32}} \rule{0pt}{4ex}\\

$A_1 \times A_1$ & $\frac{[42]_q [32]_q}{[2]_q [4]_q}$ & \rule{0pt}{4ex}\\
$A_2$ & $\frac{[42]_q [32]_q}{[2]_q [6]_q}$ & \rule{0pt}{4ex} \\
$I_2(5)$ & $\frac{[40]_q [32]_q}{[2]_q [10]_q }$ & \rule{0pt}{4ex} \\
$A_1 \times A_2, I_2(5) \times A_1, A_3, H_3$ & $\frac{[32]_q}{[2]_q}$ & \rule{0pt}{4ex} \\
\bottomrule
\end{tabular}
\caption{Sieving polynomials and divisors in type $H_4$}
\end{table}

\begin{table}[!h]
\centering
\renewcommand{\arraystretch}{1.3}
\begin{tabular}{@{} l l c @{}}
\toprule
\textbf{$\lambda$} & \textbf{$\mu_\lambda(q)$} & \textbf{Divisors of $h+2$} \\
\midrule

$A_1$ & $\frac{[20]_q [18]_q [14]_q}{[2]_q  [4]_q [6]_q}$ & \multirow{6}{*}{\parbox[c]{3cm}{\centering
1,2,7,14}} \rule{0pt}{4ex}\\

$A_1 \times A_1$ & $\frac{[18]_q [14]_q}{[2]_q [2]_q}$ & \rule{0pt}{4ex}\\
$A_2$ & $\frac{[18]_q [14]_q}{[2]_q [6]_q}$ & \rule{0pt}{4ex} \\
$B_2$ & $\frac{[16]_q [14]_q}{[2]_q [4]_q }$ & \rule{0pt}{4ex} \\
$B_3, A_2 \times A_1$ & $\frac{[14]_q}{[2]_q}$ & \rule{0pt}{4ex} \\
\bottomrule
\end{tabular}
\caption{Sieving polynomials and divisors in type $F_4$}
\end{table}

\begin{table}[!h]
\centering
\renewcommand{\arraystretch}{1.3}
\begin{tabular}{@{} l l c @{}}
\toprule
\textbf{$\lambda$} & \textbf{$\mu_\lambda(q)$} & \textbf{Divisors of $h+2$} \\
\midrule

$A_1$ & $\frac{[21]_q [20]_q [18]_q [17]_q [14]_q}{[2]_q [3]_q [4]_q [5]_q [6]_q}$ & \multirow{15}{*}{\parbox[c]{3cm}{\centering
1,2,7,14}} \rule{0pt}{4ex}\\

$A_1 \times A_1$ & $\frac{[20]_q [18]_q [17]_q [14]_q}{[2]_q [4]_q [6]_q}$ & \rule{0pt}{4ex}\\
$A_1^3$ & $\frac{[18]_q [17]_q [14]_q}{[2]_q [2]_q [3]_q }$ & \rule{0pt}{4ex} \\
$A_2 \times A_1$ & $\frac{[18]_q [17]_q [14]_q}{[2]_q [3]_q }$ & \rule{0pt}{4ex} \\
$A_3$ & $\frac{[17]_q [16]_q [14]_q}{[2]_q [4]_q }$ & \rule{0pt}{4ex} \\
$A_4$, $A_3 \times A_1$ & $\frac{[16]_q [14]_q}{[2]_q}$ & \rule{0pt}{4ex} \\
$A_2 \times A_2$ & $\frac{[18]_q [14]_q}{[2]_q [6]_q }$ & \rule{0pt}{4ex} \\
$A_2 \times A_1^2$ & $\frac{[17]_q [14]_q}{[2]_q}$ & \rule{0pt}{4ex} \\
$D_4$ & $\frac{[15]_q [14]_q}{[2]_q [3]_q}$ & \rule{0pt}{4ex} \\
$A_5, A_2^2 \times A_1$ & $\frac{[14]_q}{[2]_q}$ & \rule{0pt}{4ex} \\
$D_5, A_4 \times A_1$ & $\frac{[14]_q}{[1]_q}$ & \rule{0pt}{4ex} \\

\bottomrule
\end{tabular}
\caption{Sieving polynomials and divisors in type $E_6$}
\end{table}

\begin{table}[!h]
\centering
\renewcommand{\arraystretch}{1.3}
\begin{tabular}{@{} l l c @{}}
\toprule
\textbf{$\lambda$} & \textbf{$\mu_\lambda(q)$} & \textbf{Divisors of $h+2$} \\
\midrule

$A_1$ & $\frac{[32]_q [30]_q [28]_q [26]_q [24]_q [20]_q}{[2]_q [4]_q [6]_q [8]_q [10]_q [6]_q}$ & \multirow{20}{*}{\parbox[c]{3cm}{\centering
1,2,4,5,10,20}} \rule{0pt}{4ex}\\

$A_1 \times A_1$ & $\frac{[30]_q [28]_q [26]_q [24]_q [20]_q}{[2]_q [4]_q [6]_q [8]_q [2]_q}$ & \rule{0pt}{4ex}\\
*$A_2 \times A_1$ & $\frac{[27]_q [26]_q [24]_q [20]_q}{[2]_q [3]_q [2]_q [4]_q}$ & \rule{0pt}{4ex} \\
$(A_1^3)'$ & $\frac{[28]_q [26]_q [24]_q [20]_q}{[2]_q [6]_q [8]_q [12]_q}$ & \rule{0pt}{4ex} \\
$(A_1^3)''$ & $\frac{[28]_q [26]_q [24]_q [20]_q}{[2]_q [2]_q [4]_q [6]_q}$ & \rule{0pt}{4ex} \\
$A_3$ & $\frac{[26]_q [24]_q [24]_q [20]_q}{[2]_q [2]_q [4]_q [6]_q}$ & \rule{0pt}{4ex} \\
$(A_3 \times A_1)', A_1^4$ & $\frac{[26]_q [24]_q [20]_q}{[2]_q [4]_q [6]_q}$ & \rule{0pt}{4ex} \\
$(A_3 \times A_1)'', A_2 \times A_1^2$ & $\frac{[26]_q [24]_q [20]_q}{[2]_q [2]_q [2]_q}$ & \rule{0pt}{4ex} \\
$A_2 \times A_2$, *$A_4$ & $\frac{[26]_q [24]_q [20]_q}{[2]_q [2]_q [6]_q}$ & \rule{0pt}{4ex} \\
$D_4$ & $\frac{[24]_q [22]_q [20]_q}{[2]_q [4]_q [6]_q}$ & \rule{0pt}{4ex} \\
$A_5', A_2 \times A_1^3$ & $\frac{[24]_q [20]_q}{[2]_q [6]_q}$ & \rule{0pt}{4ex} \\
$A_5'', D_5$ & $\frac{[22]_q [20]_q}{[2]_q [2]_q}$ & \rule{0pt}{4ex} \\
$D_4 \times A_1$ & $\frac{[22]_q [20]_q}{[2]_q [4]_q}$ & \rule{0pt}{4ex} \\
$A_2 \times A_3, A_2^2 \times A_1, A_3 \times A_1^2$ & $\frac{[24]_q [20]_q}{[2]_q [2]_q}$ & \rule{0pt}{4ex} \\
$E_6, D_6, D_6 \times A_1, A_6,  A_1 \times A_5, A_4 \times A_4, A_1 \times A_2 \times A_3$ & $\frac{[20]_q}{[2]_q}$ & \rule{0pt}{4ex} \\

\bottomrule
\end{tabular}
\caption{Sieving polynomials and divisors in type $E_7$}
\end{table}

\begin{table}[!h]
\centering
\renewcommand{\arraystretch}{1.3}
\begin{tabular}{@{} l l c @{}}
\toprule
\textbf{$\lambda$} & \textbf{$\mu_\lambda(q)$} & \textbf{Divisors of $h+2$} \\
\midrule

$A_1$ & $\frac{[54]_q [50]_q [48]_q [44]_q [42]_q [38]_q [32]_q}{[2]_q [6]_q [8]_q [10]_q [12]_q [14]_q [18]_q}$ & \multirow{20}{*}{\parbox[c]{3cm}{\centering
1,2,4,8,16,32}} \rule{0pt}{4ex}\\
$A_1 \times A_1$ & $\frac{[50]_q [48]_q [44]_q [42]_q [38]_q [32]_q}{[2]_q [4]_q [6]_q [8]_q [10]_q [12]_q}$ & \rule{0pt}{4ex}\\
$A_1^3$ & $\frac{[48]_q [44]_q [42]_q [38]_q [32]_q}{[2]_q [2]_q [6]_q [8]_q [12]_q}$ & \rule{0pt}{4ex} \\
$A_3$ & $\frac{[48]_q [44]_q [42]_q [38]_q [32]_q}{[2]_q [4]_q [6]_q [8]_q [10]_q}$ & \rule{0pt}{4ex} \\
$A_1^4$ & $\frac{[44]_q [42]_q [38]_q [32]_q}{[2]_q [4]_q [6]_q [8]_q}$ & \rule{0pt}{4ex} \\
$A_2 \times A_1^2$ & $\frac{[44]_q [42]_q [38]_q [32]_q}{[2]_q [2]_q [4]_q [6]_q}$ & \rule{0pt}{4ex} \\
$A_3 \times A_1$ & $\frac{[42]_q [40]_q [38]_q [32]_q}{[2]_q [2]_q [4]_q [6]_q}$ & \rule{0pt}{4ex} \\
$D_4$ & $\frac{[42]_q [40]_q [38]_q [32]_q}{[2]_q [6]_q [8]_q [12]_q}$ & \rule{0pt}{4ex} \\
$D_4 \times A_1, D_5$ & $\frac{[38]_q [36]_q [32]_q}{[2]_q [4]_q [6]_q}$ & \rule{0pt}{4ex} \\
$A_5$ & $\frac{[38]_q [36]_q [32]_q}{[2]_q [2]_q [6]_q}$ & \rule{0pt}{4ex} \\
$A_2^2 \times A_1, A_2 \times A_1^3$ & $\frac{[42]_q [38]_q [32]_q}{[2]_q [2]_q [6]_q}$ & \rule{0pt}{4ex} \\
$A_3 \times A_1^2, A_3 \times A_2$ & $\frac{[42]_q [38]_q [32]_q}{[2]_q [2]_q [4]_q}$ & \rule{0pt}{4ex} \\
$D_6$ & $\frac{[34]_q [32]_q}{[2]_q [4]_q}$ & \rule{0pt}{4ex} \\
$E_6, D_4 \times A_2$ & $\frac{[36]_q [32]_q}{[2]_q [6]_q}$ & \rule{0pt}{4ex} \\
$A_3 \times A_2 \times A_1, A_4 \times A_2, A_4 \times A_1^2$ & $\frac{[38]_q [32]_q}{[2]_q [2]_q}$ & \rule{0pt}{4ex} \\
$A_5 \times A_1, D_5 \times A_1, A_6$ & $\frac{[36]_q [32]_q}{[2]_q [2]_q}$ & \rule{0pt}{4ex} \\
$A_2^2 \times A_1^2, A_3^2$ & $\frac{[38]_q [32]_q}{[2]_q [4]_q}$ & \rule{0pt}{4ex} \\
$E_7, D_7, A_4 \times A_2 \times A_1, A_4 \times A_3, A_6 \times A_1, D_5 \times A_2, A_7, E_6 \times A_1$ & $\frac{[32]_q}{[2]_q}$ & \rule{0pt}{4ex} \\

\bottomrule
\end{tabular}
\caption{Sieving polynomials and divisors in type $E_8$}
\end{table}

\end{document}